\documentclass[11pt,reqno]{amsart}
\usepackage[a4paper,margin=27mm]{geometry}
\usepackage[T1]{fontenc}
\usepackage{lmodern,amsmath,amssymb,amsthm,microtype,booktabs}
\usepackage{amsfonts}
\usepackage{amstext}
\usepackage{mathrsfs}
\usepackage{cite}
\usepackage[colorlinks=true,allcolors=blue]{hyperref}
\numberwithin{equation}{section}
\newtheorem{theorem}{Theorem}[section]

\newtheorem{corollary}[theorem]{Corollary}
\newtheorem{lemma}[theorem]{Lemma}
\theoremstyle{definition}
\newtheorem{example}[theorem]{Example}
\theoremstyle{remark}

\newcommand{\C}{\mathbb C}
\newcommand{\PP}{\mathbb P}
\newcommand{\OO}{\mathscr O}
\newcommand{\GG}{E^{\mathrm{GG}}}
\newcommand{\Nbar}{\overline N}
\newcommand{\JJ}{\mathcal J}
\newcommand{\II}{\mathcal I}
\newcommand{\FF}{\mathcal F}
\DeclareMathOperator{\Sym}{Sym}
\DeclareMathOperator{\Bs}{Bs}
\DeclareMathOperator{\ord}{ord}
\title[Uniform Algebraic Degeneracy of Entire Curves]{Uniform Algebraic Degeneracy of Entire Curves into Complex Projective Varieties}
\makeatletter
\renewcommand{\section}{\@startsection{section}{1}{0pt}{.7\linespacing\@plus\linespacing}{12pt}{\normalfont\scshape\centering}}
\renewcommand{\subsection}{\@startsection{subsection}{2}{0pt}{16pt plus 4pt minus 2pt}{7pt}{\normalfont\bfseries}}
\makeatother

\begin{document}
	
	\author[X.-J. Dong* and Y.-D. Lv]{Xianjing Dong* and Yanda Lv}
	
	\address{School of Mathematical Sciences \\ Qufu Normal University \\ Qufu, Jining, Shandong, 273165, P. R. China}
	
	\email{xjdong05@126.com}
	
	\address{School of Mathematical Sciences \\ Qufu Normal University \\ Qufu, Jining, Shandong, 273165, P. R. China}
	
	\email{lvyanda2000@163.com}
	
	\subjclass[2020]{32H30; 32Q45; 14F06.}
	\keywords{Nevanlinna theory; GGL conjecture;  entire curves;  logarithmic jet differentials; algebraic degeneracy.}
	\thanks{*Corresponding author}
	
	\hypersetup{pdftitle={Uniform Algebraic Degeneracy of Entire Curves into Complex Projective Varieties},pdfkeywords={Nevanlinna theory, logarithmic jet differentials, second main theorem, boundary multiplicities}}
	\date{}
	\begin{abstract}
		
		This paper is devoted to  the   algebraic degeneracy  of entire curves  intersecting a  simple-normal-crossing  divisor on smooth  projective varieties over $\mathbb C,$ within the logarithmic Green-Griffiths-Lang  program.
		For negatively‑twisted global logarithmic jet differentials, we introduce an intrinsic weighted conormal filtration along each positive-dimensional closed boundary stratum, with weights prescribed by the contact orders along  boundary components.
		Applying  the maximal Harder-Narasimhan slopes of logarithmic cotangent bundles and the normal degrees of the strata, we define an integer vanishing threshold which can improve the pole bound obtained from the slope-based ceiling.
		Under the stable-base-locus hypothesis for the logarithmic Semple tower, we  then establish an improved second main theorem whose counting coefficients depend on the  boundary stratum,  actual contact orders and the ratio of source weight to twisting exponent.
		These coefficients are independent of the powers  used to realize the stable base locus.
		We further derive a uniform algebraic‑degeneracy criterion allowing positive corrected pole coefficients, formulated via a finite system of stratum‑wise inequalities together with the ampleness of a residual divisor class. This yields conditional orbifold Green-Griffiths-Lang type statements.
		Refinements arising  from conormal‑bundle decompositions and sharp local pole estimates at surface boundary crossings are established. 
		Our examples exhibit    strict improvement over naive  slope-based  bounds
		and illustrate  that full-weight poles can persist at boundary  crossings, and  furnish   concrete instances  
		satisfying   the stable‑base‑locus hypothesis.
	\end{abstract}
	\maketitle
	\thispagestyle{empty}
	\setcounter{page}{1}
	
	\section{Introduction}\label{sec:intro}
	
	\subsection{Motivation}

	Negatively twisted jet differentials are a basic tool in the study of algebraic degeneracy and logarithmic hyperbolicity of entire curves $f:\C\to X$ intersecting a  simple-normal-crossing divisor $D,$ where $X$ is a smooth  projective variety over $\mathbb C.$  
	For an ample line bundle $L$ over $X,$ the logarithmic vanishing theorem \cite{Dem,DL,CDDR} states that 
	every  global section $$Q\in H^0\bigl(X,E_{k,m}^{\mathrm{GG}}\Omega_X^1(\log D)\otimes L^{-1}\bigr),\ \ \ \ k,m\ge1$$ 
	satisfies  $Q(j_kf)=0$ for every nonconstant entire curve $f$ omitting $D.$
	When a curve $f$ meets $D$ and $Q(j_kf)\not\equiv0,$ comparison of the zeros and poles of this evaluation yields a second main theorem in Nevanlinna theory.
	The counting coefficient is controlled by the actual poles of the evaluated differential, which the usual estimate bounds by its source weight.
	We study the reduction in this pole bound forced by the geometry of the boundary.

	When $L=H^a$ with  $a\ge1$ for an ample line bundle $H$ over $X,$   a nonzero evaluation gives a counting coefficient $m/a$ relative to the characteristic $T_f(r,H).$
	Huynh-Vu-Xie \cite[Theorem 3.1]{HVX} and Brotbek-Deng \cite[Lemma 4.7 and Corollary 4.8]{BD} formulate this principle using effective jet differentials and a stable base locus on the logarithmic Semple tower, respectively. 
	The analytic argument retains the actual pole divisor before applying  source-weight bound.
	Our question is  how much negative twisting  reduces this divisor uniformly over all global sections.
	
	A uniform improvement  to the counting coefficient requires vanishing that grows linearly with  the twisting exponent.
	We obtain this growth from a weighted conormal filtration over the boundary strata.
	The weights record contact orders and the graded factors involve  logarithmic cotangent bundle of a stratum and its normal line bundles.
	Their slopes give both an integer threshold and a simpler closed form bound, for arbitrary jet order and target dimension.
	At crossings, the bounds from individual components and from their closed intersection provide complementary constraints.
	
	The estimates apply to every section supplied by a construction which satisfies  the stable base locus hypothesis.
	Concrete families satisfying this hypothesis are available: Xie's ample-cotangent constructions yield the first-order case \cite[Theorem 1.1]{Xie}, while Brotbek--Deng establish the logarithmic higher-order condition directly \cite[Corollary 4.5]{BD}.
	Further first-order examples are provided by the almost ample logarithmic cotangent bundles of Brotbek--Deng \cite[Theorem A and Definition 2.1]{BDpos}.
	Deng's higher-order constructions \cite[Theorem C and Corollary D]{DengDT} yield further examples via the natural inclusion of ordinary invariant jet differentials into logarithmic invariant jet differentials, which preserves evaluations away from the boundary.
	The resulting coefficient improvement survives passage to powers and optimization of the same collection of weight-to-twist ratios.
	Vanishing of all contact coefficients excludes entire curves with the prescribed ramification; see Corollary~\ref{cor:ramification}.
	More generally, Theorem~\ref{thm:degeneracy} absorbs positive pole coefficients into boundary contacts and confines all such curves to one fixed proper algebraic subset.
	Corollary~\ref{cor:ggl} states the resulting conditional orbifold Green-Griffiths-Lang conclusion and specifies the remaining jet-geometric hypothesis.
	
	The general multiplicity and the analytic principles used in the argument are established.
	Logarithmic jet bundles and Semple towers are developed in \cite{DL} and logarithmic residue loci occur in \cite[Section 2.3]{BDpos}.
	Normal and tangential filtrations for logarithmic differential quotients appear in Campana-P\u aun \cite[proof of Proposition 4.7]{CP} and Cadorel \cite[Section 6]{Cad}.
	Petrakiev \cite[Lemma 2.1 and Corollary 2.2]{Pet} expresses the general conormal multiplicity principle through the first nonzero normal term of a section.
	Our construction specifies the contact-weighted filtration and its finite source-weight pieces for logarithmic jets over SNC (simple normal crossing) strata.
	Orbifold jet algebras \cite[Definition 2.5]{CDDR} impose pole conditions through prescribed ramification indices; here the bounds follow from the negative twist of an arbitrary global section.
	The slope inequalities are those of Ramanan--Ramanathan \cite[Theorem 3.18 and Section 4]{RR}, and the passage from actual poles to counting functions is the established argument in \cite{Siu,HVX,BD}.
	
	Our comparison with \cite[Corollary 4.8]{BD} fixes both the stable base locus hypothesis and the weight-to-twist data.
	The conclusion is a uniform geometric improvement  to that coefficient.
	Specialized plane-curve estimates \cite{HHMX,HHMXconics,CHLX,HWX} also depend on explicit section constructions and may give stronger numerical bounds for particular pairs.
	Some of these constructions use restricted jet loci; the additional geometry needed to treat them is discussed in Section~\ref{sec:proofs}.
	
	\subsection{Main results}
	
	Throughout the paper, $X$ is a smooth irreducible projective variety over $\C$ of dimension $n\ge2,$ and $D=\sum_{i=1}^tD_i$ is a reduced 
	simple-normal-crossing divisor with smooth irreducible components on $X.$
	Let us  consider entire curves $f:\C\to X.$
	Let $\pi:\JJ_k=J_k(X,\log D)\to X$ denote  the affine logarithmic jet space, whose 
	polynomial functions of source weight $m$ form $\GG_{k,m}\Omega_X^1(\log D),$ where an $h$th derivative has source weight $h.$
	We denote the invariant subsheaf by $E_{k,m}\Omega_X^1(\log D).$
	Its sections satisfy $$Q\bigl(j_k(\gamma\circ\varphi)\bigr)=\varphi'(0)^mQ(j_k\gamma)$$ for every invertible change of source parameter $\varphi$ fixing the origin.
	In particular, all estimates for Green-Griffiths differentials below apply to invariant differentials.
	In ideal-membership statements, $Q$ denotes the corresponding fiberwise polynomial section of $\pi^*L^{-1}$ on $\JJ_k,$ homogeneous of source weight $m.$
	We write $x_+=\max\{x,0\}$ and use the convention $\ord_a(0)=+\infty$ at every source point $a.$
	Orders of meromorphic line-bundle sections are computed in a holomorphic frame, so they do not depend on that frame.
	For such a section $s,$ put $\operatorname{pole}_a s=(-\ord_a s)_+$ when $s\not\equiv0,$ and $\operatorname{pole}_a(0)=0.$
	The notation $\mathscr O_T$ denotes the structure sheaf of a variety $T;$ when $F\simeq\PP^1,$ the notation $\mathscr O_F(b)$ denotes the line bundle of degree $b.$
	
	Fix a nonempty index set $J$ and a connected component $Y$ of $D_J=\bigcap_{i\in J}D_i.$
	Write $$Y^\circ=Y\setminus\bigcup_{i\notin J}D_i, \ \ \ \ B_Y=\sum_{i\notin J}D_i\Big{|}_Y,\ \ \ \ N_i=\OO_X(D_i)\big{|}_Y\ \ (i\in J).$$
	When $\dim Y>0,$ the divisor $B_Y$ is reduced and simple normal crossing on $Y.$
	For an ample polarization $A_Y$ on $Y$ of dimension $d_Y>0$ and a nonzero torsion-free coherent sheaf $\mathscr E,$ put $$\deg_{A_Y}\mathscr E=c_1(\mathscr E)A_Y^{d_Y-1},\ \ \ \ \mu_{A_Y}(\mathscr E)=\frac{\deg_{A_Y}\mathscr E}{\operatorname{rank}\mathscr E}.$$
	The maximal Harder--Narasimhan slope is denoted by $\mu_{\max,A_Y}(\mathscr E).$
	For the zero sheaf,  we set $\mu_{\max,A_Y}(0)=-\infty.$
	For a line bundle $L$ over $X,$ we use the numerical data $$\ell_Y=\deg_{A_Y}(L|_Y),\ \ \ \ \rho_Y=\mu_{\max,A_Y}\bigl(\Omega_Y^1(\log B_Y)\bigr),\ \ \ \ \nu_{Y,i}=\deg_{A_Y}N_i.$$
	For a contact vector $\mathbf p=(p_i)_{i\in J}$ of positive integers, set \begin{equation}\label{eq:mu}\mu_Y(\mathbf p)=\max\bigl\{0,\rho_Y,-\nu_{Y,i}/p_i:i\in J\bigr\}.\end{equation}
	The polarization is fixed on each closed stratum, independently of any entire curve or contact order.
	
	Over a positive-dimensional closed stratum $Y,$ we let $Z_{Y,k}$ be the pure normal jet locus: it retains only the logarithmic jet variables normal to the branches indexed by $J.$
	It is canonically $Y\times\mathbb A^{|J|k};$ write $\II_Y$ for its ideal and $\xi$ for its normal jet parameters.
	Let $w_i$ denote a local defining function of $D_i.$
	For a contact vector $\mathbf p,$ write $\FF_{Y,\mathbf p}^{r}$ for the local monomial filtration assigning degree $p_i$ to a defining function of $D_i,$ degree one to every tangential jet variable, and degree zero to each retained normal jet variable.
	The intrinsic construction and its coordinate invariance are established in Section~\ref{sec:prelim}.
	For $\mathbf p=\mathbf1=(1)_{i\in J},$ this is the ordinary power $\II_Y^r.$
	
	For a positive-dimensional closed stratum $Y$ with $\ell_Y>0,$ define the integer threshold \begin{equation}\label{eq:sigma}\sigma_Y(L,\mathbf p)=\min\left\{j+\sum_{i\in J}p_i b_i:\ j\in\mathbb Z_{\ge0},\ \mathbf b\in\mathbb Z_{\ge0}^{J},\ \ j\rho_Y-\sum_{i\in J}b_i\nu_{Y,i}\ge\ell_Y\right\}.\end{equation}
	If the set is empty, then  denote  this minimum   by $+\infty$.
	This is a necessary numerical threshold for a graded term to carry a section; its attainment is not an existence assertion.
	
	Throughout the paper, the notation $``\|"$ indicates that Nevanlinna's inequalities hold for all $r>1$ outside a set of finite Lebesgue measure. 
	The implied constants and exceptional sets may depend on the fixed geometric data, section, and curve, but are chosen once along the curve and are not varied at individual boundary contacts.

	The first result controls the actual poles of every negatively twisted global jet differential.

	\begin{theorem}[Weighted boundary multiplicity]\label{thm:mult}
		Let $Y$ be a positive-dimensional closed stratum  with  $\ell_Y>0$ as above,  and let $\mathbf p$ have positive integer entries. When the threshold is finite, 
		every global section $$Q\in H^0\bigl(X,\GG_{k,m}\Omega_X^1(\log D)\otimes L^{-1}\bigr), \ \ \ \     k,m\ge1$$
		belongs to $\FF_{Y,\mathbf p}^{\,\sigma_Y(L,\mathbf p)}\otimes\pi^*L^{-1}$. 
		If the threshold is infinite, then $Q=0.$
		For every holomorphic arc $\gamma$ centered at $Y^\circ,$ not contained in $D,$ with $\ord_0(w_i\circ\gamma)\ge p_i$ for $i\in J,$ one has \begin{equation}\label{eq:order}\ord_0Q(j_k\gamma)\ge\sigma_Y(L,\mathbf p)-m.\end{equation}
		In particular, if $\mu_Y(\mathbf p)>0,$ then $$\sigma_Y(L,\mathbf p)\ge\left\lceil\frac{\ell_Y}{\mu_Y(\mathbf p)}\right\rceil,$$ and the pole order is at most $$\left(m-\left\lceil\frac{\ell_Y}{\mu_Y(\mathbf p)}\right\rceil\right)_+.$$
		If $\mu_Y(\mathbf p)=0,$ then every such global section is zero.
	\end{theorem}
	
	The integer threshold records a discrete optimization and can be strictly larger than the simpler slope ceiling.
	For a fixed source weight,  one may retain the additional constraint $$\sigma_Y^{[m]}(L,\mathbf p)=\min\left\{j+\sum_i p_i b_i:0\le j\le m,\ j\in\mathbb Z,\ \mathbf b\in\mathbb Z_{\ge0}^{J},\ j\rho_Y-\sum_i b_i\nu_{Y,i}\ge\ell_Y\right\},$$ and again assigning the value $+\infty$ to an empty feasible set.
	One has $\sigma_Y^{[m]}(L,\mathbf p)\ge\sigma_Y(L,\mathbf p),$ and every occurrence of $\sigma_Y$ in Theorem~\ref{thm:mult} can be replaced by this stronger finite-weight threshold.
	In particular, an empty constrained feasible set forces the corresponding weight-$m$ space of global sections to vanish.
	Nonemptiness of either feasible set is only a necessary numerical condition and does not imply the existence of a section or attainment of its predicted vanishing order.
	Example~\ref{ex:ruled} realizes a strict difference between the integer threshold and the slope ceiling on a logarithmic surface with nonzero negatively twisted sections.
	The proof of Theorem~\ref{thm:mult} is deferred to Section~\ref{sec:proofs}.
	
	Fix an ample line bundle $H$ over $X$  with a smooth Hermitian metric $h_H$ of positive curvature $\omega_H=c_1(H, h_H).$
	We use $d^c=\frac{\sqrt{-1}}{2\pi}(\bar\partial-\partial),$ so that $dd^c=\frac{\sqrt{-1}}{\pi}\partial\bar\partial$ and $c_1(L, h_L)=-dd^c\log\|e\|_{h_L}$ for a local holomorphic frame $e$ of a metrized line bundle $L.$ 
	For $r\ge1,$ define $$T_f(r,H)=\int_1^r\frac{dt}{t}\int_{|z|<t}f^*\omega_H.$$
	For any  line bundle $L$ over  $X$ with a fixed smooth metric $h_L,$   $T_f(r,L)$ is defined by the same formula with $c_1(L,h_L).$
	We use additive notation for line-bundle and rational-divisor classes and extend the characteristic to rational divisors by rational linearity.
	A change of smooth metric changes the characteristic by $O(1).$
	For an effective real divisor $A=\sum_z\lambda_z[z]$ on $\C,$ put $$n_A(t)=\sum_{|z|<t}\lambda_z,\ \ \ \ N(r,A)=\int_1^r\frac{n_A(t)}{t}\,dt.$$
	This convention includes the multiplicity at the origin and differs from other standard lower-limit conventions only by a bounded term.
	We use $\log^+x=\max\{\log x,0\}$ for $x>0$ and set $\log^+0=0.$
	For a nonconstant entire curve $f:\C\to X$ not contained in $D,$ write $$S_f(r)=O\big(\log^+T_f(r,H)+\log r\big),\ \ \ \ \Nbar_f(r,D)=N\bigl(r,(f^*D)_{\mathrm{red}}\bigr).$$
	
	Let $\pi_{0,k}:X_k(D)\to X$ denote the logarithmic Demailly--Semple tower and $X_k(D)^{\mathrm{sing}}$ its nonregular-jet locus.
	For a line bundle $M$ on any projective variety, put $\mathbf B(M)=\bigcap_{q\ge1}\Bs|M^q|.$
	For each positive-dimensional closed stratum,  choose a polarization as above and set $$h_Y=\deg_{A_Y}(H|_Y)>0,\ \ \ \ \theta_Y(\mathbf p)=\frac{h_Y}{\mu_Y(\mathbf p)}\ \   \  (\mu_Y(\mathbf p)>0).$$
	The positivity needed for these ratios follows from the hypothesis of Theorem~\ref{thm:base}.
	For a component $D_i,$ abbreviate this function by $\theta_i(p_i).$
	The data on $D_i$ use its full residual boundary $(D-D_i)|_{D_i}.$
	
	Whenever these ratios are defined, for $c>0,$ a connected exact stratum $Y^\circ\subset D_J,$ and a contact vector $\mathbf p,$ set \begin{eqnarray}\label{eq:coefficient}b_Y(\mathbf p)&=&\left(\min\left\{c,\ |J|c-\sum_{i\in J}\theta_i(p_i),\ c-\theta_Y(\mathbf p)\right\}\right)_+,\ \ \ \ \dim Y>0,\nonumber\\b_Y(\mathbf p)&=&\left(\min\left\{c,\ |J|c-\sum_{i\in J}\theta_i(p_i)\right\}\right)_+,\ \ \ \ \dim Y=0.\end{eqnarray}
	At a point $z\in f^{-1}(D),$ let $Y(z)^\circ$ be the exact stratum containing $f(z)$ and let $\mathbf p(z)$ be the vector of contact orders.
	Define the effective real divisor $$\mathcal B_f=\sum_{z\in f^{-1}D}b_{Y(z)}(\mathbf p(z))[z].$$
	The coefficients lie between zero and $c,$ so its counting function is well-defined.
	
	The second one converts this uniform algebraic estimate into a second main theorem under the usual stable-base-locus hypothesis.
	
	\begin{theorem}[Uniform second main theorem]\label{thm:base}
		Suppose that  $H$ is ample and there are positive integers $k,m,a$ such that
		\begin{equation}\label{eq:base}\mathbf B\big(\mathscr O_{X_k(D)}(m)\otimes\pi_{0,k}^*H^{-a}\big)\subset X_k(D)^{\mathrm{sing}}\cup\pi_{0,k}^{-1}(D).\end{equation}
		Then, the quantities $\mu_Y(\mathbf p)$ are positive for all positive-dimensional strata and finite contact vectors.
		With $c=m/a$ and the coefficients \eqref{eq:coefficient}, every nonconstant entire curve $f:\C\to X$ not contained in $D$ satisfies \begin{equation}\label{eq:uniform}T_f(r,H)\le N(r,\mathcal B_f)+S_f(r) \ \ \    \|. 
		\end{equation}
		The coefficients require no condition on the residues of selected sections and are independent of the power used to realize the stable base locus.
	\end{theorem}
	
	Every coefficient in \eqref{eq:coefficient} is at most the coefficient $c$ in \cite[Corollary 4.8]{BD}.
	On an exact stratum of positive dimension it is strictly smaller than $c.$
	At zero-dimensional crossings equality can occur.
	The quantities $\theta_i(p_i)$ and $\theta_Y(\mathbf p)$ are determined by $H,$ the polarized boundary strata and the contact orders, while the corrected coefficient $b_Y(\mathbf p)$ also depends on $c=m/a.$
	The estimate imposes no further condition on the sections supplied by \eqref{eq:base}.
	
	For a smooth irreducible boundary, put $h=\deg_A(H|_D),$ $\rho=\mu_{\max,A}(\Omega_D^1),$ $\nu=\deg_A N_{D/X},$ and $\mu(p)=\max\{\rho,-\nu/p,0\}.$
	Theorem~\ref{thm:base} gives, in particular, $$T_f(r,H)\le\left(\frac ma-\frac h{\mu(1)}\right)_+\Nbar_f(r,D)+S_f(r) \ \ \    \|.$$
	The positive subtraction survives every power used to realize the stable base locus.
	Corollary~\ref{cor:smooth} gives the finer estimate using actual contact orders, and Corollary~\ref{cor:ramification} identifies when the resulting pole cancellation excludes nonconstant entire curves.
	
	Fix integers $r_i\ge1.$ It is said  that a curve has minimum contacts $\mathbf r$,  if every intersection with $D_i$ has order at least $r_i.$
	Omission of a component satisfies this condition.
	For an exact stratum indexed by $J,$ we write $\mathbf r_J=(r_i)_{i\in J}.$
	
	The third one allows a fixed algebraic exceptional set and gives a finite-contact criterion for uniform algebraic degeneracy.
	
	\begin{theorem}[Uniform algebraic degeneracy]\label{thm:degeneracy}
		Let $Z\subsetneq X$ be a closed algebraic subset, possibly empty. 
		Suppose that   \begin{equation}\label{eq:baseZ}\mathbf B\bigl(\mathscr O_{X_k(D)}(m)\otimes\pi_{0,k}^*H^{-a}\bigr)\subset X_k(D)^{\mathrm{sing}}\cup\pi_{0,k}^{-1}(D\cup Z).\end{equation}
		for positive integers $k,m,a.$
		With $c=m/a,$ the ratios in \eqref{eq:coefficient} are well-defined, and \eqref{eq:uniform} holds for every nonconstant entire curve not contained in $D\cup Z.$
		Suppose  that  there are rational numbers $\lambda_i\ge0$ such that  
		\begin{equation}\label{eq:absorption}b_Y(\mathbf r_J)\le\sum_{i\in J}\lambda_i r_i,\ \ \ \ H-\sum_{i=1}^t\lambda_iD_i\ \text{is ample} \end{equation}
		for every nonempty exact stratum $Y^\circ$ indexed by $J.$
		Then,  every nonconstant entire curve not contained in $D$ with minimum contacts $\mathbf r$ has image contained in $Z.$
		In particular, the exceptional set is independent of the curve, and  there exist  no such curves if $Z$ is empty.
	\end{theorem}
	
	There are only finitely many stratum inequalities in \eqref{eq:absorption}.
	The theorem allows $b_Y(\mathbf r_J)>0$ and therefore goes beyond the zero-pole test of Corollary~\ref{cor:ramification}.
	The passage from pole absorption to degeneracy uses the established vanishing mechanism for negatively twisted invariant jets; the additional numerical information is the corrected coefficient $b_Y$ and its dependence on the exact stratum and contacts.
	For comparison at fixed $k,m,a,$ the full-weight pole estimate requires the stronger inequalities $c\le\sum_{i\in J}\lambda_i r_i.$
	Neither comparison supplies the base locus condition itself.
	
	\subsection{Further refinements}
	
	The weighted filtration bounds slopes by the largest slope of its normal and tangential factors.
	A nonsplit conormal bundle can give a stronger unweighted estimate.
	The following refinement applies the conormal multiplicity principle to the full extension bundle.
	
	\begin{theorem}[A conormal fiber and branch partitions]\label{prop:fiber}
		Let $Y$ be positive-dimensional with $\ell_Y>0,$ and set $$E_\xi=\left.(\II_Y/\II_Y^2)\right|_{Y\times\{\xi\}},\ \ \ \ \eta_\xi=\max\{0,\mu_{\max,A_Y}(E_\xi)\}.$$
		For any chosen normal jet parameter $\xi_0$ with $\eta_{\xi_0}>0,$ every section in Theorem~\ref{thm:mult} belongs to $$\II_Y^{\lceil\ell_Y/\eta_{\xi_0}\rceil}\otimes\pi^*L^{-1}.$$
		If $\eta_{\xi_0}=0,$ then every such section is zero.
		In particular, let $\mathcal P=(G_1,\ldots,G_k)$ be a partition $J=G_1\sqcup\cdots\sqcup G_k,$ allowing empty groups.
		Let $X_G$ be the component of $D_{J\setminus G}$ containing $Y,$ with $X_G=X$ when $G=J,$ and put $$R_G=\sum_{i\notin J}D_i\Big|_{X_G},\ \ \ \ F_G=\Omega_{X_G}^1(\log R_G)\Big|_Y.$$
		For $G=\varnothing,$ one puts $F_G=\Omega_Y^1(\log B_Y).$
		There is a parameter $\xi_0$ with $E_{\xi_0}\simeq\bigoplus_{h=1}^kF_{G_h},$ so the preceding conclusions hold with 
		$$\eta_{\xi_0}=\eta_{\mathcal P}=\max\big\{0,\mu_{\max,A_Y}(F_{G_h}): 1\le h\le k\big\}\le\mu_Y(\mathbf1).$$
	\end{theorem}
	
	For a component $D_i,$ abbreviate $\II_{D_i}$ by $\II_i.$
	On a surface, the simultaneous constraints from two boundary components give the following local estimate.
	
	\begin{theorem}[Exact bound from two component ideals]\label{prop:cross}
		Let $k,m\ge1,$  $s_i,s_j\in\mathbb Z_{\ge0},$ and suppose that a weight-$m$ logarithmic $k$-jet differential on a surface belongs locally to $\II_i^{s_i}\cap\II_j^{s_j}$ at $D_i\cap D_j.$
		Here,  ideal membership is required in the fiberwise polynomial jet algebra over an adapted neighborhood of the target crossing.
		Let $\gamma$ be an arc through this crossing, not contained in $D,$ with finite contacts $p_i,p_j\in\mathbb Z_{\ge1}.$
		Put $$b_0=(s_i-m)_+,\ \ \ \ d_0=(s_j-m)_+,\ \ \ \ \delta=(s_i+s_j-m)_+,$$ and $$C=p_i b_0+p_jd_0+\min(p_i,p_j)(\delta-b_0-d_0).$$
		Then, $\delta\ge b_0+d_0$ and the evaluation has pole order at most $(m-C)_+.$
		This is the best local bound obtainable from these two monomial ideal-power conditions and source weight alone.
	\end{theorem}
	
	\subsection{Structure of the article}
	
	The statements in Section~\ref{sec:intro} are proved in Section~\ref{sec:proofs}.
	Section~\ref{sec:prelim} supplies the slope inequality, the logarithmic-jet filtration, the established actual-pole estimate, and a pole-absorption lemma.
	Section~\ref{sec:proofs} also contains uniform algebraic degeneracy and conditional orbifold Green-Griffiths-Lang consequences, explicit surface and hypersurface computations, and examples testing the integer threshold and the crossing bounds.
	
	\section{Preliminaries}\label{sec:prelim}
	
	\subsection{Slopes and vanishing}
	
	To include applications to boundary strata, let $X$ in this subsection denote any smooth irreducible projective variety of dimension $d\ge1$ over $\C,$ equipped with an ample line bundle $A.$
	Slopes are computed with respect to $A,$ and the subscript $A$ is omitted when the polarization is clear.
	A coherent sheaf $\mathscr E$ is torsion-free if multiplication by every nonzero element of $\mathscr O_{X,y}$ is injective on $\mathscr E_y$ for every $y\in X.$
	Its rank is the dimension of its generic fiber over the function field of $X.$
	For a torsion-free sheaf $\mathscr E$ of positive rank $r,$ define its determinant, degree, and slope by $$\det\mathscr E=\left(\bigwedge\nolimits^r\mathscr E\right)^{**},\ \ \ \ \deg_A\mathscr E=c_1(\det\mathscr E)c_1(A)^{d-1},\ \ \ \ \mu_A(\mathscr E)=\frac{\deg_A\mathscr E}{r},$$ where the double dual is taken with respect to $\mathscr O_X.$
	Here, $c_1$ is the first Chern class and the product is taken on $X.$
	Smoothness ensures that $\det\mathscr E$ is a line bundle; on a curve, $\deg_A\mathscr E$ is its usual degree.
	
	A subsheaf $\mathscr F\subset\mathscr E$ is saturated if $\mathscr E/\mathscr F$ is torsion-free.
	Its saturation is the inverse image in $\mathscr E$ of the torsion subsheaf of $\mathscr E/\mathscr F.$
	Write $\mathscr F^{\mathrm{sat}}$ for the saturation of a nonzero subsheaf $\mathscr F.$
	The codimension-one lengths of $\mathscr F^{\mathrm{sat}}/\mathscr F$ define an effective divisor $\Delta,$ and $$\deg_A\mathscr F^{\mathrm{sat}}-\deg_A\mathscr F=c_1(\mathscr O_X(\Delta))c_1(A)^{d-1}\ge0.$$
	Thus, saturation preserves rank and cannot decrease degree, so upper bounds for slopes may be tested on saturated subsheaves.
	A nonzero torsion-free sheaf $\mathscr E$ is slope semistable if $$\mu_A(\mathscr F)\le\mu_A(\mathscr E)$$ for every subsheaf $\mathscr F$ with $0<\operatorname{rank}\mathscr F<\operatorname{rank}\mathscr E.$
	For every nonzero torsion-free sheaf $\mathscr E,$ the Harder-Narasimhan filtration is the unique filtration by saturated subsheaves $$0=\mathscr E_0\subsetneq\mathscr E_1\subsetneq\cdots\subsetneq\mathscr E_s=\mathscr E$$ such that the quotients $\mathscr G_i=\mathscr E_i/\mathscr E_{i-1}$ are torsion-free and slope semistable, with $$\mu_A(\mathscr G_1)>\mu_A(\mathscr G_2)>\cdots>\mu_A(\mathscr G_s).$$
	The maximal slope is $$\mu_{\max,A}(\mathscr E)=\mu_A(\mathscr G_1)=\sup_{0\ne\mathscr F\subset\mathscr E}\mu_A(\mathscr F),$$ where the supremum is taken over all nonzero coherent subsheaves, including those of full rank.
	As in the introduction, $\mu_{\max,A}(0)=-\infty;$ the slope itself is used only for positive-rank sheaves.
	In particular, semistability is equivalent to $\mu_{\max,A}(\mathscr E)=\mu_A(\mathscr E).$
	
	A torsion-free sheaf is locally free outside a closed subset of codimension at least two.
	If two torsion-free sheaves agree on such an open set, they have the same rank and first Chern class, and each saturated subsheaf of either sheaf extends to a saturated subsheaf of the other with the same rank and degree.
	Their maximal slopes therefore coincide, and semistability is preserved.
	We shall use this observation to compute the slopes of filtration factors on a common open set where they are locally free.
	
	For an exact sequence of positive-rank torsion-free sheaves $0\to\mathscr A\to\mathscr B\to\mathscr C\to0,$ taking determinants gives $$\det\mathscr B\simeq\det\mathscr A\otimes\det\mathscr C,\ \ \ \ \deg_A\mathscr B=\deg_A\mathscr A+\deg_A\mathscr C.$$
	The determinant isomorphism is obtained on the open set where all three sheaves are locally free and extends across its complement of codimension at least two.
	Similarly, let $\mathscr A,\mathscr B$ have positive ranks $a,b,$ and let $\mathscr T$ be the quotient of $\mathscr A\otimes\mathscr B$ by its torsion subsheaf.
	The determinant formula for a tensor product of vector bundles extends from such an open set to give $\det\mathscr T\simeq(\det\mathscr A)^{\otimes b}\otimes(\det\mathscr B)^{\otimes a}.$
	Consequently,  we get 
	$$\deg_A\mathscr T=b\deg_A\mathscr A+a\deg_A\mathscr B,\ \ \ \ \mu_A(\mathscr T)=\mu_A(\mathscr A)+\mu_A(\mathscr B).$$
	
	\begin{lemma}[Slope inequalities]\label{lem:slope}
		Let $E,F$ be nonzero vector bundles over $X,$ and $P$  a line bundle over $X.$ For $j\ge1,$ we have 
		\begin{eqnarray*}\mu_{\max}(E\otimes F)&\le&\mu_{\max}(E)+\mu_{\max}(F),\\ \mu_{\max}(\Sym^jE)&\le&j\mu_{\max}(E),\\ \mu_{\max}(E\otimes P)&=&\mu_{\max}(E)+\deg_A P,\\ \mu_{\max}(E\oplus F)&=&\max\{\mu_{\max}(E),\mu_{\max}(F)\}.\end{eqnarray*}
		For an exact sequence of torsion-free sheaves $0\to\mathscr E'\to\mathscr E\to\mathscr E''\to0,$ one has 
		$$\mu_{\max}(\mathscr E)\le\max\big\{\mu_{\max}(\mathscr E'), \   \mu_{\max}(\mathscr E'')\big\}.$$
		A torsion-free sheaf of negative maximal slope has no nonzero global section.
		We use the convention $\Sym^0E=\mathscr O_X.$
	\end{lemma}
	
	\begin{proof}
		We begin with the assertion for an exact sequence of torsion-free sheaves.
		The extension inequality is immediate if $\mathscr E=0,$ so assume $\mathscr E\ne0.$
		For a nonzero subsheaf $\mathscr F\subset\mathscr E,$ put $\mathscr F'=\mathscr F\cap\mathscr E'$ and let $\mathscr F''$ be its image in $\mathscr E''.$
		These are torsion-free, and $$0\longrightarrow\mathscr F'\longrightarrow\mathscr F\longrightarrow\mathscr F''\longrightarrow0$$ is exact.
		When both terms are nonzero, with ranks $r',r'',$  the additivity of rank and degree gives $$\mu_A(\mathscr F)=\frac{r'}{r'+r''}\mu_A(\mathscr F')+\frac{r''}{r'+r''}\mu_A(\mathscr F'').$$
		If one of these terms is zero, the slope of $\mathscr F$ equals the slope of the other term.
		Each of those slopes is bounded by the maximal slope of its ambient sheaf.
		Taking the supremum over $\mathscr F$ proves the extension inequality.
		In characteristic zero, the quotient by torsion of the tensor product of two slope-semistable torsion-free sheaves is slope semistable; see \cite[Theorem 3.18 and its higher-dimensional extension in Section 4]{RR}.
		Write $E_i,F_q$ for the terms of the two Harder--Narasimhan filtrations and $\mathcal G_i=E_i/E_{i-1},$ $\mathcal H_q=F_q/F_{q-1}$ for their factors.
		Choose an open set $U\subset X$ whose complement has codimension at least two and on which all these terms and factors are locally free.
		First filter $(E\otimes F)|_U$ by $W_i=(E_i\otimes F)|_U.$
		The quotient map $$p_i:W_i\longrightarrow(\mathcal G_i\otimes F)|_U$$ has kernel $W_{i-1}.$
		If $0=F_0\subset\cdots\subset F_t=F,$ then we refine this step by $$W_{i,q}=p_i^{-1}\bigl((\mathcal G_i\otimes F_q)|_U\bigr),\ \ \ \ 0\le q\le t.$$
		Then $$W_{i,0}=W_{i-1},\ \ \ \ W_{i,t}=W_i,\ \ \ \ W_{i,q}/W_{i,q-1}\simeq(\mathcal G_i\otimes\mathcal H_q)|_U.$$
		Concatenating these refinements gives a filtration of $(E\otimes F)|_U$ by subbundles.
		Let $\eta$ be the generic point of $X.$
		For each term, define $$\widetilde W_{i,q}=(E\otimes F)\cap(W_{i,q})_\eta\subset(E\otimes F)_\eta,$$ where the intersection is taken inside the sheaf of rational sections of $E\otimes F.$
		This intersection is coherent, as is seen locally by clearing denominators in linear equations defining $(W_{i,q})_\eta.$
		The quotient $(E\otimes F)/\widetilde W_{i,q}$ embeds in the sheaf of rational sections with fiber $(E\otimes F)_\eta/(W_{i,q})_\eta,$ so it is torsion-free.
		Thus, $\widetilde W_{i,q}$ is saturated and restricts to $W_{i,q}$ on $U,$ since $W_{i,q}$ is a subbundle there.
		The extensions remain nested, and their successive quotients $$\mathcal Q_{i,q}=\widetilde W_{i,q}/\widetilde W_{i,q-1},\ \ \ \ 1\le q\le t,$$ are torsion-free since  they are subsheaves of $(E\otimes F)/\widetilde W_{i,q-1}.$
		They satisfy $$\mathcal Q_{i,q}|_U\simeq(\mathcal G_i\otimes\mathcal H_q)|_U.$$
		Let $\mathcal T_{i,q}$ be the quotient of $\mathcal G_i\otimes\mathcal H_q$ by its torsion subsheaf.
		By the cited characteristic-zero tensor product theorem, $\mathcal T_{i,q}$ is slope semistable, and the tensor product slope formula gives $$\mu_A(\mathcal T_{i,q})=\mu_A(\mathcal G_i)+\mu_A(\mathcal H_q).$$
		Since $\mathcal Q_{i,q}$ and $\mathcal T_{i,q}$ agree on $U,$ invariance under changes in codimension at least two gives $$\mu_{\max}(\mathcal Q_{i,q})=\mu_{\max}(\mathcal T_{i,q})=\mu_A(\mathcal G_i)+\mu_A(\mathcal H_q)\le\mu_{\max}(E)+\mu_{\max}(F).$$
		Applying the extension inequality successively proves the tensor product bound.
		The averaging operator $$\frac{1}{j!}\sum_{\sigma\in\mathfrak S_j}\sigma,$$ where the symmetric group $\mathfrak S_j$ permutes the tensor factors, exhibits $\Sym^jE$ as a direct summand of $E^{\otimes j}.$
		The symmetric power bound follows from this inclusion and the iterated tensor product bound.
		Tensoring by $P$ gives a bijection between subsheaves of $E$ and subsheaves of $E\otimes P.$
		It preserves rank and adds $\deg_A P$ to every slope, proving the twisting formula.
		The extension inequality gives the upper bound for $E\oplus F,$ while the inclusions of its two summands give the reverse bound.
		Finally, a nonzero section of a torsion-free sheaf $\mathscr E$ is nonzero at the generic point and hence gives an injection $\mathscr O_X\hookrightarrow\mathscr E.$
		Since $\mu_A(\mathscr O_X)=0,$ this forces $\mu_{\max}(\mathscr E)\ge0$ and proves the vanishing assertion.
	\end{proof}
	
	Negative average slope alone does not imply vanishing: $E=\mathscr O_{\PP^1}(1)\oplus\mathscr O_{\PP^1}(-3)$ has $\mu(E)=-1$ and $\mu_{\max}(E)=1,$ with $H^0(\PP^1,E)\ne0.$
	For a finite filtration $0=\mathcal K_0\subset\cdots\subset\mathcal K_N=\mathcal K$ with nonzero torsion-free quotients $\mathcal Q_v=\mathcal K_v/\mathcal K_{v-1},$ repeated application of Lemma~\ref{lem:slope} gives $$\mu_{\max}(\mathcal K)\le\max_{1\le v\le N}\mu_{\max}(\mathcal Q_v).$$
	If all these maximal slopes are negative, then the maximum is negative because the filtration is finite, and $H^0(X,\mathcal K)=0.$
	Equivalently, the left-exact sequences $$0\longrightarrow H^0(X,\mathcal K_{v-1})\longrightarrow H^0(X,\mathcal K_v)\longrightarrow H^0(X,\mathcal Q_v)$$ give the same vanishing by induction.
	
	We apply these facts to a positive-dimensional stratum $Y,$ with polarization $A_Y$ and numerical data $\rho_Y,\nu_{Y,i},\ell_Y$ as in the introduction.
	For nonnegative integers $b_i,\alpha_h,$ set $$\mathcal G_{\mathbf b,\boldsymbol\alpha}=\left(\bigotimes_{i\in J}N_i^{-b_i}\right)\otimes\left(\bigotimes_{h=1}^k\Sym^{\alpha_h}\Omega_Y^1(\log B_Y)\right)\otimes L^{-1}|_Y.$$
	The tensor product, symmetric power, and twisting formulas give $$\mu_{\max,A_Y}(\mathcal G_{\mathbf b,\boldsymbol\alpha})\le\left(\sum_{h=1}^k\alpha_h\right)\rho_Y-\sum_{i\in J}b_i\nu_{Y,i}-\ell_Y.$$
	Here, $N_i^{-b_i}=(N_i^*)^{\otimes b_i},$ so its degree is $-b_i\nu_{Y,i},$ and the twist $L^{-1}|_Y$ contributes $-\ell_Y.$
	So, the strict negativity of the right-hand side forces $H^0(Y,\mathcal G_{\mathbf b,\boldsymbol\alpha})=0.$
	The integer $j=\sum_h\alpha_h$ counts the factors from $\Omega_Y^1(\log B_Y),$ whereas their total source weight is $\sum_hh\alpha_h.$
	For a factor arising from a weight-$m$ jet differential, this gives $j\le\sum_hh\alpha_h\le m.$
	This estimate underlies \eqref{eq:sigma} and its refinement at fixed source weight.
	Thus, a nonzero section of $\mathcal G_{\mathbf b,\boldsymbol\alpha}$ requires $$\left(\sum_h\alpha_h\right)\rho_Y-\sum_i b_i\nu_{Y,i}\ge\ell_Y;$$ this necessary condition alone does not imply the existence of a section.
	
	\subsection{Logarithmic jet coordinates and filtrations}
	
	Let $Y$ denote a positive-dimensional closed stratum indexed by $J,$ with residual boundary $B_Y.$
	The logarithmic cotangent bundle $\Omega_Y^1(\log B_Y)$ has a local frame consisting of $dz_a/z_a$ in boundary directions and ordinary differentials in the remaining directions.
	For $i\in J,$ the normal line bundle of $D_i$ restricted to $Y$ is $N_i=\OO_X(D_i)|_Y,$ and its dual is $N_i^*=\OO_X(-D_i)|_Y.$
	Choose adapted coordinates $(w_i,z_a)$ near a point of $Y,$ where $w_i=0$ defines $D_i$ for $i\in J$ and $1\le a\le d_Y.$
	After reordering the $z_a,$ the residual boundary branches are $z_1=0,\ldots,z_q=0,$ where $0\le q\le d_Y.$
	For the branches indexed by $J,$ use the logarithmic jet coordinates $$\xi_{i,h}=(\log w_i)^{(h)},\ \ \ \ i\in J,\ \ \ \ 1\le h\le k.$$
	These agree with the indicated derivatives away from the boundary and extend as independent fiber coordinates over $w_i=0.$
	Write $\xi_i=(\xi_{i,1},\ldots,\xi_{i,k})$ and denote the remaining jet variables by $\eta_{a,h},$ with $1\le a\le d_Y$ and $1\le h\le k.$
	Explicitly, $$\eta_{a,h}=\begin{cases}(\log z_a)^{(h)},&1\le a\le q,\\ z_a^{(h)},&q<a\le d_Y.\end{cases}$$
	Thus, the base coordinates together with $\xi_{i,h},\eta_{a,h}$ give a local trivialization of the affine logarithmic jet space, with $nk$ independent fiber coordinates.
	Set $$\II_Y=(w_i,\eta_{a,h}:i\in J,\ 1\le a\le d_Y,\ 1\le h\le k).$$
	Locally, its zero locus $Z_{Y,k}$ lies over $Y,$ with $\eta_{a,h}=0$ and unrestricted coordinates $\xi_{i,h}.$
	The ideals are considered on neighborhoods of this locus in the logarithmic jet space.
	For positive integer contact weights $\mathbf p=(p_i)_{i\in J},$ define in each adapted chart \begin{equation}\label{eq:filtration}\FF_{Y,\mathbf p}^{r}=\left(w^{\mathbf b}\eta^{\boldsymbol\alpha}:\sum_{i\in J}p_i b_i+|\boldsymbol\alpha|\ge r\right),\ \ \ \ r\in\mathbb Z_{\ge0},\end{equation} where $w^{\mathbf b}=\prod_{i\in J}w_i^{b_i},$ $\eta^{\boldsymbol\alpha}=\prod_{a,h}\eta_{a,h}^{\alpha_{a,h}},$ $|\boldsymbol\alpha|=\sum_{a,h}\alpha_{a,h},$ and all exponents are nonnegative integers.
	The filtration degrees of $w_i,\xi_{i,h},\eta_{a,h}$ are $p_i,0,1,$ respectively.
	Under the source rescaling $t\mapsto\lambda t,$ the base coordinates remain fixed and $\xi_{i,h},\eta_{a,h}$ are multiplied by $\lambda^h;$ the source weights of $w_i,\xi_{i,h},\eta_{a,h}$ are therefore $0,h,h.$
	The filtration is decreasing and multiplicative: $$\FF_{Y,\mathbf p}^0=\OO_{\JJ_k},\ \ \ \ \FF_{Y,\mathbf p}^1=\II_Y,\ \ \ \ \FF_{Y,\mathbf p}^{r+1}\subset\FF_{Y,\mathbf p}^r,\ \ \ \ \FF_{Y,\mathbf p}^{r}\FF_{Y,\mathbf p}^{s}\subset\FF_{Y,\mathbf p}^{r+s}.$$
	Equivalently, with $\mathfrak a_i=(w_i),$ one has $$\FF_{Y,\mathbf p}^{r}=\sum_{\sum_i p_i b_i+j\ge r}\left(\prod_{i\in J}\mathfrak a_i^{b_i}\right)\II_Y^j.$$
	Indeed, every monomial in $\II_Y^j$ has filtration degree at least $j,$ and every generator $w^{\mathbf b}\eta^{\boldsymbol\alpha}$ of \eqref{eq:filtration} occurs in the displayed sum with $j=|\boldsymbol\alpha|.$
	Writing $p_{\max}=\max_{i\in J}p_i,$ one also has $$\II_Y^r\subset\FF_{Y,\mathbf p}^{r}\subset\II_Y^{\lceil r/p_{\max}\rceil}.$$
	The first inclusion uses $p_i\ge1,$ and the second one follows because the filtration degree of a monomial is at most $p_{\max}$ times its total degree in the generators of $\II_Y.$
	Assigning filtration degree one to every $\eta_{a,h}$ is compatible with arbitrary positive contact weights.
	Indeed, if $p_i=1$ and $h\ge2,$ the term $w_iB_h(\xi_i)$ in the transition below has filtration degree one, so assigning degree $h$ to $\eta_{a,h}$ would not give a coordinate-independent filtration in general.
	
	\begin{lemma}\label{lem:filtration}
		The ideals $\II_Y$ define a regular closed embedding $$Z_{Y,k}\simeq Y\times\mathbb A^{|J|k}\ \hookrightarrow\ \JJ_k,$$ with ideal understood on a neighborhood of this closed locus.
		If $\pi_Y:Z_{Y,k}\to Y$ is the projection, then there is an exact sequence \begin{equation}\label{eq:conormal}0\longrightarrow\bigoplus_{i\in J}\pi_Y^*N_i^*\longrightarrow\II_Y/\II_Y^2\longrightarrow\bigoplus_{h=1}^k\pi_Y^*\Omega_Y^1(\log B_Y)[h]\longrightarrow0.\end{equation}
		The symbol $[h]$ places this copy of $\Omega_Y^1(\log B_Y)$ in source weight $h,$ while the classes of $w_i$ have source weight zero.
		The filtration \eqref{eq:filtration} is independent of the adapted coordinates.
		For fixed $r,m\in\mathbb Z_{\ge0},$ the homogeneous component $$\left(\pi_{Y*}\bigl(\FF_{Y,\mathbf p}^{r}/\FF_{Y,\mathbf p}^{r+1}\bigr)\right)_m$$ has a finite filtration whose quotients are finite direct sums of $$\left(\bigotimes_{i\in J}N_i^{-b_i}\right)\otimes\left(\bigotimes_{h=1}^k\Sym^{\alpha_h}\Omega_Y^1(\log B_Y)\right).$$
		The indices $b_i,\alpha_h$ are nonnegative integers satisfying $$\sum_i p_i b_i+\sum_h\alpha_h=r,\ \ \ \ \sum_hh\alpha_h\le m.$$
		The multiplicity is the dimension of the space of polynomials of source weight $m-\sum_hh\alpha_h$ in the variables $\xi_{i,h}.$
		For $\mathbf p=\mathbf1,$ one has $\FF_{Y,\mathbf1}^r=\II_Y^r.$
	\end{lemma}
	
	\begin{proof}
		Under a change of adapted coordinates,  we have $w_i'=u_iw_i$ with $u_i$ a unit.
		Let $B_h$ be the complete Bell polynomial determined by $$w_i^{(h)}=w_iB_h(\xi_{i,1},\ldots,\xi_{i,h}).$$
		For example $$B_1=\xi_{i,1},\ \ \ \ B_2=\xi_{i,1}^2+\xi_{i,2},\ \ \ \ B_3=\xi_{i,1}^3+3\xi_{i,1}\xi_{i,2}+\xi_{i,3}.$$
		Every $B_h$ is homogeneous of source weight $h.$
		Let $D_t$ denote total differentiation from jet order $h$ to order $h+1$ for $0\le h<k;$ all expressions below have resulting jet order at most $k.$
		The identities $$D_tw_i=w_i\xi_{i,1},\ \ \ \ D_t\xi_{i,h}=\xi_{i,h+1},\ \ \ \ D_t\eta_{a,h}=\eta_{a,h+1}\ \ (h<k)$$ show that membership in the corresponding ideals generated by $w_i,\eta_{a,h}$ at successive orders, and in their squares, is preserved under differentiation.
		For the assertion about squares, the Leibniz rule gives $$D_t(g_1g_2)=(D_tg_1)g_2+g_1D_tg_2,$$ and both summands belong to the square of the ideal at the next jet order whenever $g_1,g_2$ belong to the ideal at the preceding order.
		Moreover, the total derivative of every function of the base coordinates belongs to this ideal, since $D_tz_a=z_a\eta_{a,1}$ for $a\le q$ and $D_tz_a=\eta_{a,1}$ otherwise.
		Consequently, $$\xi'_{i,h}=\xi_{i,h}+D_t^h\log u_i\equiv\xi_{i,h}\pmod{\II_Y}.$$
		The derivatives of $\log u_i$ are regular expressions in the unit $u_i$ and its derivatives and do not require a chosen logarithm.
		Write $\vartheta_a=dz_a/z_a$ for $a\le q$ and $\vartheta_a=dz_a$ otherwise.
		Under a change of adapted coordinates, these forms satisfy $$\vartheta'_a=\sum_b\widetilde A_{ab}(w,z)\vartheta_b+\sum_{i\in J}\widetilde C_{ai}(w,z)\,dw_i,$$ with regular coefficients.
		For a residual boundary branch this follows from $z'_a=v_az_a$ and $d\log z'_a=d\log z_a+d\log v_a,$ after a possible permutation of the branches.
		Restricting to $Y$ gives the transition matrix $A=\widetilde A|_Y$ of $\Omega_Y^1(\log B_Y).$
		Evaluation on a 1-jet gives $$\eta'_{a,1}=\sum_b\widetilde A_{ab}\eta_{b,1}+\sum_{i\in J}\widetilde C_{ai}w_i\xi_{i,1}.$$
		When this identity is differentiated, derivatives of its coefficients belong to $\II_Y$ and are multiplied by another element of $\II_Y.$
		These products vanish modulo $\II_Y^2.$
		Using $D_t(w_iB_h(\xi_i))=w_iB_{h+1}(\xi_i),$ induction on $h$, it therefore gives \begin{eqnarray*}w_i'&\equiv&u_i|_Y\,w_i,\\ \eta'_{a,h}&\equiv&\sum_b A_{ab}\eta_{b,h}+\sum_{i\in J}C_{ai}w_iB_h(\xi_i)\pmod{\II_Y^2},\end{eqnarray*} where $C_{ai}=\widetilde C_{ai}|_Y.$
		In particular, all transformed generators of $\II_Y$ belong to $\II_Y,$ and the inverse change of coordinates gives the reverse inclusion.
		The local ideals therefore agree on overlaps, and the coordinates $\xi_{i,h}$ have identity transition maps on $Z_{Y,k}.$
		The generators $w_i,\eta_{a,h}$ form part of a local coordinate system, so they form a regular sequence.
		The resulting regular embedding has codimension $|J|+kd_Y.$
		The classes of $w_i$ therefore span the conormal subbundle $\bigoplus_i\pi_Y^*N_i^*,$ and the quotient is given by the tangential generators of each source weight $h.$
		This proves \eqref{eq:conormal} without assuming a splitting.
		For $\mathbf p=\mathbf1,$ the monomial definition gives $\FF_{Y,\mathbf1}^r=\II_Y^r$ for every $r\ge0.$
		For $\mathbf p=\mathbf1,$ regularity of the embedding also gives the canonical isomorphism $$\bigoplus_{r\ge0}\II_Y^r/\II_Y^{r+1}\simeq\Sym_{\OO_{Z_{Y,k}}}^{\bullet}(\II_Y/\II_Y^2).$$
		The exact sequence \eqref{eq:conormal} induces filtrations on these symmetric powers; a global splitting of the conormal bundle is not required.
		To verify that the entire filtration is well-defined, apply the chain rule before reducing modulo $\II_Y^2.$
		Every term in a transformed tangential jet contains at least one tangential jet factor or at least one $w_i,$ since that jet vanishes on $Z_{Y,k}.$
		It therefore has filtration degree at least one.
		The function $w_i'=u_iw_i$ has filtration degree at least $p_i,$ and each transformed coordinate $\xi'_{i,h}$ has filtration degree at least zero.
		Substitution consequently takes every generator of \eqref{eq:filtration} into the same filtration level.
		The inverse change of coordinates has the same properties, so the resulting filtration is independent of the adapted coordinates.
		At every point $\zeta\in Z_{Y,k},$ the inclusion $\FF_{Y,\mathbf p}^r\subset\II_Y^{\lceil r/p_{\max}\rceil}$ and the Krull intersection theorem give $$\bigcap_{r\ge0}(\FF_{Y,\mathbf p}^r)_\zeta=0.$$
		Indeed, $\II_{Y,\zeta}$ is contained in the maximal ideal of the Noetherian local ring $\mathscr O_{\JJ_k,\zeta}.$
		The same assertion holds after tensoring by a line bundle, so a nonzero germ at $\zeta$ has a finite largest filtration order and a nonzero image in the corresponding graded quotient.
		By $\II_Y\FF_{Y,\mathbf p}^{r}\subset\FF_{Y,\mathbf p}^{r+1},$ each associated graded quotient is naturally a sheaf on $Z_{Y,k}.$
		We work throughout with the algebra of functions that are polynomial in the jet variables.
		The jet space and ideal sheaves in this argument are considered in the algebraic category.
		The homogeneous jet bundles on $X$ and the graded bundles with both filtration degree and source weight fixed on $Y$ are coherent; GAGA on these projective bases identifies their algebraic and holomorphic global sections \cite{SerreGAGA}.
		The affine direct-image argument below is therefore applied to quasi-coherent algebraic sheaves.
		Write $$\mathcal A=\bigoplus_{r\ge0}\FF_{Y,\mathbf p}^{r}/\FF_{Y,\mathbf p}^{r+1}.$$
		Over an adapted affine open set $U\subset Y,$ its algebra of sections is $$\mathcal A(\pi_Y^{-1}(U))\simeq\mathscr O_Y(U)[\xi_{i,h}][W_i,V_{a,h}],$$ where $W_i$ and $V_{a,h}$ are the initial classes of $w_i$ and $\eta_{a,h},$ of filtration degrees $p_i$ and one and source weights zero and $h,$ respectively.
		There are no relations among these initial classes, since $w_i,\eta_{a,h}$ are coordinate generators.
		This calculation may also be made in adapted \'etale coordinates, since flat base change preserves the monomial ideals and their associated graded quotients.
		The preceding transition formulas imply $$W_i'=u_i|_Y\,W_i,\ \ \ \ V'_{a,h}=\sum_bA_{ab}V_{b,h}+\sum_{i\in J: \ p_i=1}C_{ai}W_iB_h(\xi_i).$$
		The first formula uses the exact equality $w_i'=u_iw_i,$ rather than only its reduction modulo $\II_Y^2,$ and the degree-zero initial class of $u_i$ is $u_i|_Y.$
		Thus, it remains valid for arbitrary $p_i,$ including $p_i>1.$
		For the second formula, the terms discarded modulo $\II_Y^2$ have filtration degree at least two, because $\II_Y^2\subset\FF_{Y,\mathbf p}^{2}.$
		Only the terms with $p_i=1$ appear in the second formula, because $V_{a,h}$ has filtration degree one.
		Let $\mathcal K\subset\mathcal A$ be the homogeneous ideal generated locally by $W_i$ for $i\in J.$
		The transition formula for $W_i$ shows that this ideal is independent of the adapted coordinates.
		Both $\mathcal K$ and its powers are homogeneous for the filtration degree and the source weight.
		For fixed $r,$ the degree-$r$ parts $(\mathcal K^s)_r$ define a finite decreasing filtration of $\mathcal A_r,$ with $(\mathcal K^s)_r=0$ for $s>r.$
		In the associated graded algebra for this second filtration, the terms containing $W_iB_h(\xi_i)$ disappear from the transformation of $V_{a,h}.$
		Thus, $W_i$ transforms as a frame of $\pi_Y^*N_i^*,$ while the variables $V_{a,h}$ for each fixed $h$ transform as a frame of $\pi_Y^*\Omega_Y^1(\log B_Y).$
		The latter transformations do not mix different values of $h.$
		For a local monomial, put $\alpha_h=\sum_a\alpha_{a,h};$ this is the number of tangential factors of derivative order $h.$
		The coordinate transformations within that order yield $\Sym^{\alpha_h}\Omega_Y^1(\log B_Y).$
		It follows that $(\mathcal K^s)_r/(\mathcal K^{s+1})_r$ is the direct sum of the bundles $$\pi_Y^*\left[\left(\bigotimes_{i\in J}N_i^{-b_i}\right)\otimes\left(\bigotimes_{h=1}^k\Sym^{\alpha_h}\Omega_Y^1(\log B_Y)\right)\right]$$ indexed by $$\sum_i b_i=s,\ \ \ \ \sum_i p_i b_i+\sum_h\alpha_h=r.$$
		The tensor generators in such a summand have source weight $\sum_hh\alpha_h.$
		Since the coordinates $\xi_{i,h}$ have identity transition maps on $Z_{Y,k},$ its weight-$m$ direct image is the underlying bundle on $Y$ tensored with $$\C[\xi_{i,h}]_{\,m-\sum_hh\alpha_h},$$ where the subscript denotes the homogeneous part for source weights $\operatorname{wt}(\xi_{i,h})=h.$
		This is the homogeneous polynomial space on the affine factor of $Z_{Y,k}\simeq Y\times\mathbb A^{|J|k};$ it is zero for negative weight and finite-dimensional for each nonnegative weight.
		Moreover, fixed $r$ bounds all $b_i,\alpha_h,$ so only finitely many summands occur.
		These monomials form a local basis over $\mathscr O_Y,$ so each fixed-$(r,m)$ homogeneous component is a vector bundle of finite rank.
		The projection $\pi_Y$ is affine, so its direct image is exact on these quasi-coherent sheaves.
		All maps preserve source weight, and passage to the homogeneous component of a fixed source weight is also exact.
		This proves the finite filtration assertion after direct image.
	\end{proof}
	
	On $Z_{Y,k},$ the restriction of $\pi^*L^{-1}$ is $\pi_Y^*(L^{-1}|_Y).$
	This line bundle has source weight zero, and the projection formula gives \begin{eqnarray*}&&\left(\pi_{Y*}\left[(\FF_{Y,\mathbf p}^r/\FF_{Y,\mathbf p}^{r+1})\otimes\pi_Y^*(L^{-1}|_Y)\right]\right)_m\\ &\simeq&\left(\pi_{Y*}(\FF_{Y,\mathbf p}^r/\FF_{Y,\mathbf p}^{r+1})\right)_m\otimes L^{-1}|_Y.\end{eqnarray*}
	Its filtration factors are the bundles $\mathcal G_{\mathbf b,\boldsymbol\alpha}$ from the preceding subsection, with the multiplicities of Lemma~\ref{lem:filtration}.
	If all these factors have negative maximal slope, Lemma~\ref{lem:slope} gives the vanishing of the global sections of this homogeneous component.
	For a global weight-$m$ differential $Q$ lying in $\FF_{Y,\mathbf p}^r\otimes\pi^*L^{-1}$ along $Z_{Y,k},$ its graded image is a global section over the whole closed stratum $Y.$
	If the preceding vanishing applies, that image is zero, and left exactness of global sections gives $Q\in\FF_{Y,\mathbf p}^{r+1}\otimes\pi^*L^{-1}$ along $Z_{Y,k}.$
	Starting from $\FF_{Y,\mathbf p}^0=\mathscr O_{\JJ_k},$ iteration gives $Q\in\FF_{Y,\mathbf p}^N\otimes\pi^*L^{-1}$ whenever the twisted weight-$m$ graded components in degrees $0,\ldots,N-1$ have no global sections.
	If this vanishing holds in every degree, then local separatedness forces the germ of $Q$ at each point of $Z_{Y,k}$ to be zero.
	The jet space is irreducible, being locally an affine-space bundle over the irreducible variety $X,$ so the identity theorem gives $Q=0.$
	
	\subsection{Evaluation of   logarithmic jet differentials along entire curves}
	
	Throughout this subsection, $f:\C\to X$ is a nonconstant entire curve with $f(\C)\not\subset D.$
	We use the characteristic, counting functions, and error term defined in Section~\ref{sec:intro}.
	We evaluate jets using the fixed global coordinate $z$ on $\C,$ so $j_kf(z)$ is the jet at $t=0$ of the germ $t\mapsto f(z+t).$
	The transition laws of the Green-Griffiths bundle make the evaluation independent of the target coordinates.
	It is holomorphic outside $f^{-1}(D)$ and extends meromorphically across this discrete set by the local expressions in logarithmic derivatives.
	Reparametrization invariance will be used in Lemma~\ref{lem:absorption} to combine the evaluation with the source differential on $\PP^1.$
	When $s=Q(j_kf)\not\equiv0,$ its pole and zero divisors are $$P_Q=\sum_{a\in\C}(-\ord_a s)_+[a],\ \ \ \ R_Q=\sum_{a\in\C}(\ord_a s)_+[a].$$
	The orders are computed in local holomorphic frames of $f^*L^{-1}$ and do not depend on those frames.
	Both sums are locally finite, and the support of $P_Q$ is contained in $f^{-1}(D).$
	The following estimate retains these divisors with their actual multiplicities.
	
	\begin{lemma}\label{lem:analytic}
		For $$Q\in H^0\bigl(X,\GG_{k,m}\Omega_X^1(\log D)\otimes L^{-1}\bigr),\ \ \ \ Q(j_kf)\not\equiv0,$$ one has \begin{equation}\label{eq:analytic}T_f(r,L)+N(r,R_Q)\le N(r,P_Q)+S_f(r) \ \ \   \|.
		\end{equation}
	\end{lemma}
	
	\begin{proof}
		Use the fixed smooth Hermitian metric $h_L$ defining $T_f(r,L).$
		In a local holomorphic frame $e$ of $L,$ write $s=q(f^*e)^{-1},$ so that $$\log\|s\|=\log|q|-\log\|f^*e\|_{h_L}.$$
		With the normalization fixed in the introduction, the Poincar\'e--Lelong formula gives $$dd^c\log\|s\|=[R_Q]-[P_Q]+f^*c_1(L,h_L).$$
		The logarithmic singularities at zeros and poles are locally integrable, and their circular means are finite even when a zero or pole lies on the circle.
		Put $$M_s(r)=\int_0^{2\pi}\log\|s(re^{i\vartheta})\|\frac{d\vartheta}{2\pi}.$$
		Integrating the current identity over discs and then over the radii from $1$ to $r$ gives Jensen's formula in the form $$T_f(r,L)+N(r,R_Q)-N(r,P_Q)=M_s(r)-M_s(1).$$
		The identity may first be obtained for circles avoiding zeros and poles and then extended to all radii by continuity of the circular means.
		The term $M_s(1)$ is a finite constant; the zeros and poles at the origin are already included in the counting functions.
		Choose finitely many adapted coordinate neighborhoods $V_\alpha$ and relatively compact open subsets $U_\alpha\Subset V_\alpha$ covering $X,$ with coordinate functions that are restrictions of rational functions.
		Such charts are obtained by completing the local algebraic equations of the boundary branches to a regular system of parameters, shrinking analytically, and then using compactness to select a finite cover.
		Trivialize $L$ on each $V_\alpha;$ the coefficients of $Q,$ the metric factors, and the ordinary coordinate functions are bounded on $\overline U_\alpha.$
		Discard every index $\alpha$ for which $U_\alpha\cap f(\C)=\varnothing.$
		For each remaining coordinate function $g,$ choose $x\in U_\alpha\cap f(\C).$
		Since $g$ is regular at $x,$ it admits a representation $g=P_g/S_g$ by sections of a sufficiently large power $H^{b_g},$ with $S_g(x)\ne0.$
		Indeed, the pole divisor $F_g$ of $g$ does not contain $x,$ and $H^{b_g}(-F_g)$ is globally generated for $b_g$ sufficiently large.
		Choose a section $u$ of this line bundle nonzero at $x,$ and let $\sigma_{F_g}$ be the canonical section of $\mathscr O_X(F_g).$
		Then, $S_g=u\sigma_{F_g}$ and $P_g=gS_g$ are sections of $H^{b_g}$ with the required properties.
		Here, $gS_g$ has no divisorial poles and extends holomorphically over $X,$ since $X$ is smooth and hence normal.
		Thus, $f^*S_g\not\equiv0,$ and $g_f=g\circ f$ defines a meromorphic function on $\C.$
		Let $C_g$ be the effective divisor whose ideal sheaf is locally generated by the coefficients of $f^*P_g$ and $f^*S_g$ in a common holomorphic frame.
		If $f^*P_g\equiv0,$ then $C_g$ is the zero divisor of $f^*S_g.$
		Use the metric induced by $f^*h_H^{b_g}$ to define $$\psi_g=\frac12\log\left(\|f^*P_g\|^2+\|f^*S_g\|^2\right).$$
		After cancellation of these common zeros, the pair defines $g_f:\C\to\PP^1.$
		If $\omega_{\mathrm{FS}}$ denotes  the Fubini--Study form representing $c_1(\mathscr O_{\PP^1}(1)),$ Poincar\'e--Lelong gives $$dd^c\psi_g=g_f^*\omega_{\mathrm{FS}}+[C_g]-b_gf^*\omega_H.$$
		Here, $T_{g_f}$ is computed with respect to $\omega_{\mathrm{FS}};$ it agrees with the usual meromorphic characteristic up to $O(1).$
		Applying Jensen's formula to this identity gives $$T_{g_f}(r)+N(r,C_g)-b_gT_f(r,H)=\int_0^{2\pi}\bigl(\psi_g(re^{i\vartheta})-\psi_g(e^{i\vartheta})\bigr)\frac{d\vartheta}{2\pi}.$$
		The function $\psi_g$ is bounded above by compactness of $X,$ and its circular mean at radius one is finite, so $$T_{g_f}(r)+N(r,C_g)\le b_gT_f(r,H)+O(1).$$
		In particular, $T_{g_f}(r)\le b_gT_f(r,H)+O(1),$ as required for the logarithmic derivative estimate.
		For a boundary defining function, $g_f\not\equiv0$ follows from $f(\C)\not\subset D$ and the identity theorem.
		Each evaluated logarithmic jet coordinate is a polynomial in the quotients $g_f^{(q)}/g_f.$
		An ordinary derivative, when $g_f\not\equiv0,$ is written as $g_f^{(h)}=(g_f^{(h)}/g_f)g_f.$
		If an ordinary coordinate has $g_f\equiv0,$ all its derivatives vanish and are omitted.
		For each nonzero $g_f$ and $1\le q\le k,$ the logarithmic derivative lemma gives $$\int_0^{2\pi}\log^+\left|\frac{g_f^{(q)}}{g_f}(re^{i\vartheta})\right|\frac{d\vartheta}{2\pi}\le C_{g,q}\bigl(\log^+T_{g_f}(r)+\log r\bigr) \ \ \  \|.$$
		For a nonzero constant $g_f,$ all these derivative quotients vanish and the same estimate is immediate.
		The polynomial expression for $Q$ and the bounds on each relatively compact chart imply the pointwise estimate $$\log^+\|s(z)\|\le C+C\sum_g\sum_{q=1}^k\log^+\left|\frac{g_f^{(q)}(z)}{g_f(z)}\right|,$$ away from the discrete sets where a quotient is undefined.
		The sum ranges over the finitely many coordinate functions used in the cover whose compositions with $f$ are not identically zero.
		The extra factors $g_f$ from ordinary derivatives are bounded on $f^{-1}(U_\alpha),$ so they contribute only to the constant in this estimate.
		Taking circular means and using the rational-function estimate therefore gives $$\int_0^{2\pi}\log^+\big\|s(re^{i\vartheta})\big\|\frac{d\vartheta}{2\pi}=S_f(r) \ \ \    \|,$$
		compared with  \cite[proof of Theorem 3.1]{HVX} and \cite[Lemma 4.7]{BD}.
		Only finitely many functions and derivative orders occur, so the union of their exceptional sets has finite measure.
		The implied constant may depend on $Q$ and $f,$ but is uniform in $r$ and in the choice of a boundary contact point along this fixed curve.
		Combining this estimate with Jensen's formula and $\log\|s\|\le\log^+\|s\|$ proves \eqref{eq:analytic}.
		At a boundary contact, each logarithmic derivative of order $h$ has pole order at most $h.$
		Indeed, let $a\in f^{-1}(D),$ put $t=z-a,$ and for each boundary branch through $f(a)$ write $$w_i(f(a+t))=t^{p_i}u_i(t),\ \ \ \ p_i\ge1,\ \ \ \ u_i(0)\ne0.$$
		Then $$\frac{d^h}{dt^h}\log\bigl(w_i(f(a+t))\bigr)=(-1)^{h-1}(h-1)!\,p_i t^{-h}+\frac{d^h}{dt^h}\log u_i(t),$$ and the second term is holomorphic near $t=0.$
		The ordinary jet coordinates are holomorphic there, and at a crossing the contributions from all boundary branches in any one monomial are bounded by the sum of their derivative weights.
		Since every monomial of $Q$ has source weight $m,$ this also gives $P_Q\le m(f^*D)_{\mathrm{red}},$ and hence an inequality in terms of the reduced boundary counting function.
	\end{proof}
	
	The preceding inequality must also be used when $T_f(r,H)=O(\log r),$ in which case its error term need not be $o(T_f(r,H)).$
	The next lemma includes this case.
	
	\begin{lemma}\label{lem:absorption}
		Let $Q\in H^0(X,E_{k,M}\Omega_X^1(\log D)\otimes L^{-1})$ have positive weight $M,$ and let $R$ be an effective rational divisor supported on $D$ such that $L-R$ is ample.
		Then, there exist  no nonconstant entire curves $f$ not contained in $D$ for which $$Q(j_kf)\not\equiv0,\ \ \ \ P_Q\le f^*R.$$
	\end{lemma}
	
	\begin{proof}
		Suppose that $f$ is a nonconstant entire curve satisfying the stated conditions.
		Lemma~\ref{lem:analytic} and $P_Q\le f^*R$ give $$T_f(r,L)\le N(r,f^*R)+S_f(r) \ \ \    \|.$$
		Choose an integer $q>0$ such that $qR$ is integral.
		Scale a smooth metric on $\mathscr O_X(qR)$ so that its canonical section, whose zero divisor is $qR,$ has norm at most one on the compact variety $X.$
		The first main theorem applied to this section and divided by $q$ gives $$N(r,f^*R)\le T_f(r,R)+O(1),$$ since $f$ is not contained in the support of $R.$
		Comparison with the fixed metrics contributes only a bounded term.
		Set $B=L-R.$
		Choose a positive curvature representative $\omega_B$ of this ample rational class by dividing the curvature form of a positive integral multiple by that multiple.
		Compactness gives $\omega_B\ge\varepsilon\omega_H$ for some $\varepsilon>0,$ and hence $$\varepsilon T_f(r,H)\le T_f(r,B)+O(1)=T_f(r,L)-T_f(r,R)+O(1).$$
		Thus, 
		for some constant $C>0,$ one has $$T_f(r,H)\le C\bigl(\log^+T_f(r,H)+\log r\bigr) \ \ \     \|.$$
		The inequality  $C\log^+x\le x/2+C'$ with $x\ge0$ leads to  that $T_f(r,H)=O(\log r)$ for all $r>1$  outside a set $E$ of finite Lebesgue measure.
		For  sufficiently large $r,$ the measure of $E\cap[r,\infty)$ is less than one and so $[r,r+1]$ contains a radius outside $E.$
		Monotonicity of $T_f(r,H)$ then extends the estimate to all sufficiently large radii.
		More explicitly, for $r'\in[r,r+1]\setminus E$ and all sufficiently large $r,$ one has $$T_f(r,H)\le T_f(r',H)\le C_0\log r'\le C_0\log(r+1)=O(\log r).$$
		Choose a positive integer $a$ such that $H^a$ defines an embedding $\iota:X\hookrightarrow\PP^N.$
		The characteristic of $\iota\circ f$ is $aT_f(r,H)+O(1)=O(\log r).$
		Choose a homogeneous coordinate that does not vanish identically on the image.
		Denote this coordinate by $Z_0$ and put $g_j=(Z_j/Z_0)\circ\iota\circ f$ for $1\le j\le N.$
		The comparison for pairs of sections used in the proof of Lemma~\ref{lem:analytic} gives $$T_{g_j}(r)\le T_{\iota\circ f}(r)+O(1)=O(\log r),$$ with common zeros cancelled when defining the ratio.
		The ratios of the remaining homogeneous coordinates to this one are meromorphic functions on $\C$ of characteristic $O(\log r),$ and hence are rational by the rationality criterion; see \cite{NW}.
		Clearing denominators gives a rational map from $\PP^1$ to $\PP^N,$ and removing common factors in a local coordinate extends this map across every point.
		Its image lies in the closed subvariety $\iota(X),$ so it defines a nonconstant morphism $\bar f:\PP^1\to X$ extending $f.$
		In a fixed local frame of $L^{-1},$ let $q_z,q_\tau$ denote the evaluation coefficients in source coordinates $z,\tau.$
		Under $z=\varphi(\tau),$ reparametrization invariance gives $$q_\tau(\tau)=\varphi'(\tau)^M q_z(\varphi(\tau)),\ \ \ \ q_\tau(\tau)(d\tau)^M=q_z(z)(dz)^M.$$
		The expressions $Q(j_k\bar f)(d\tau)^M$ therefore define a nonzero meromorphic section of $$K_{\PP^1}^{\otimes M}\otimes\bar f^*L^{-1}.$$
		On the affine chart $\C\subset\PP^1,$ the source frame $(dz)^M$ has no zeros or poles, so its finite pole orders are exactly those of the original evaluation.
		Since $\bar f$ is not contained in $D,$ the pullback $\bar f^*R$ is an effective rational divisor, and the finite poles have total order at most $\deg\bar f^*R.$
		Let $\tau$ be a local parameter vanishing at infinity.
		For a local equation $w_i$ of a boundary branch through $\bar f(\infty),$ write $$w_i\circ\bar f=\tau^{p_i}u_i(\tau),\ \ \ \ p_i\ge1,\ \ \ \ u_i(0)\ne0.$$
		Then $$\frac{d^h}{d\tau^h}\log(w_i\circ\bar f)=(-1)^{h-1}(h-1)!\,p_i\tau^{-h}+\frac{d^h}{d\tau^h}\log u_i(\tau).$$
		The second term is holomorphic, so this logarithmic jet coordinate has pole order at most $h.$
		The ordinary jet coordinates, and the logarithmic coordinates for branches not passing through $\bar f(\infty),$ are holomorphic there.
		Every monomial of $Q$ has total source weight $M,$ and its coefficient is holomorphic in a local frame of $L^{-1}.$
		Consequently, $\ord_{\tau=0}q_\tau\ge-M$ in the chosen frame, including when $\bar f(\infty)$ lies on several boundary branches.
		Since $(d\tau)^M$ is a holomorphic nowhere-vanishing local frame of $K_{\PP^1}^{\otimes M}$ near infinity, this is also a bound for the pole order of the global meromorphic section at infinity.
		If $P$ is the pole divisor of this global meromorphic section on $\PP^1,$ these two bounds give $$\deg P\le\deg\bar f^*R+M.$$
		The evaluation is a nonzero holomorphic section after twisting by $\OO_{\PP^1}(P).$
		Since $\deg K_{\PP^1}=-2,$ its line bundle has degree \begin{eqnarray*}\deg\left(K_{\PP^1}^{\otimes M}\otimes\bar f^*L^{-1}\otimes\mathscr O_{\PP^1}(P)\right)&=&-2M-\deg\bar f^*L+\deg P\\ &\le&-M-\deg\bar f^*(L-R)<0.\end{eqnarray*}
		The strict inequality follows from the ampleness of $L-R$ and the nonconstancy of $\bar f.$
		This contradicts with the existence of a nonzero holomorphic section of a line bundle of negative degree.
	\end{proof}
	
	\section{Proofs of the main results and consequences}\label{sec:proofs}
	
	\subsection{Proofs of Theorems~\ref{thm:mult}, \ref{thm:base} and \ref{thm:degeneracy} and consequences}~
	
	\begin{proof}[Proof of Theorem~\ref{thm:mult}]~
		
		~
		
		Tensor a graded quotient in Lemma~\ref{lem:filtration} by $L^{-1}|_Y.$
		Writing $j=\sum_h\alpha_h,$ Lemma~\ref{lem:slope} bounds its maximal slope by $$j\rho_Y-\sum_i b_i\nu_{Y,i}-\ell_Y.$$
		For filtration degree $r<\sigma_Y(L,\mathbf p),$ this number is negative for every possible quotient, by \eqref{eq:sigma}.
		Every quotient of the finite filtration in Lemma~\ref{lem:filtration} therefore has no nonzero global section, and its consecutive exact sequences give the same vanishing for the weight-$m$ graded piece.
		Starting with $\FF_{Y,\mathbf p}^{0}=\mathscr O_{\JJ_k},$ the image of $Q$ in each consecutive graded piece is therefore zero for $r<\sigma_Y(L,\mathbf p).$
		More precisely, once $Q$ belongs to $\FF_{Y,\mathbf p}^{r}\otimes\pi^*L^{-1}$ along $Z_{Y,k},$ its graded image is a global section of the twisted weight-$m$ component just shown to vanish, so exactness gives membership in $\FF_{Y,\mathbf p}^{r+1}\otimes\pi^*L^{-1}$ there.
		Induction proves the required membership along the whole closed locus $Z_{Y,k}.$
		The indices in Lemma~\ref{lem:filtration} satisfy $j=\sum_h\alpha_h\le\sum_hh\alpha_h\le m.$
		Restricting the feasible set to $j\le m$ consequently gives the identical induction with $\sigma_Y^{[m]}$ in place of $\sigma_Y.$
		
		For the simpler slope threshold, the definition \eqref{eq:mu} gives $$j\rho_Y-\sum_i b_i\nu_{Y,i}\le\left(j+\sum_i p_i b_i\right)\mu_Y(\mathbf p).$$
		If $\mu_Y(\mathbf p)>0,$ every feasible tuple in \eqref{eq:sigma} consequently has degree at least $\lceil\ell_Y/\mu_Y(\mathbf p)\rceil.$
		If either feasible set is empty, the corresponding induction gives membership in every $\FF_{Y,\mathbf p}^{r}.$
		With $p_{\max}=\max_i p_i,$ the local monomial description gives $$\II_Y^r\subset\FF_{Y,\mathbf p}^{r}\subset\II_Y^{\lceil r/p_{\max}\rceil}.$$
		Krull intersection gives zero in each local ring along $Z_{Y,k}.$
		Thus, the section vanishes on a nonempty open set of the irreducible jet space $\JJ_k,$ and hence identically.
		The case $\mu_Y(\mathbf p)=0$ has an empty feasible set because $\ell_Y>0.$
		
		If $Q=0,$ then  the remaining order and pole bounds follow immediately from the stated conventions, so we assume $Q\ne0.$
		At an arc centered in $Y^\circ,$ all tangential jet variables are holomorphic, whereas $$\xi_{i,h}(j_k\gamma)=(-1)^{h-1}(h-1)!p_i' t^{-h}+O(1),\ \ \ \ p_i'=\ord_0(w_i\circ\gamma)\ge p_i.$$
		Let $s<+\infty$ be either threshold for which membership has been established above.
		In a fixed local holomorphic frame of $L^{-1}$ near the center of the arc, write $$Q=\sum_{\alpha,\beta}q_{\alpha\beta}(w,z)\eta^\alpha\xi^\beta,$$ where the sum is finite and the coefficients depend only on the base coordinates.
		Here, $\alpha=(\alpha_{a,h})$ and $\beta=(\beta_{i,h})$ are nonnegative multi-indices, $|\alpha|=\sum_{a,h}\alpha_{a,h},$ and source homogeneity gives $$\sum_{a,h}h\alpha_{a,h}+\sum_{i,h}h\beta_{i,h}=m$$ for every term.
		Membership along the whole locus $Z_{Y,k}$ implies the coefficient conditions $$q_{\alpha\beta}\in\left(w^{\mathbf b}:\sum_i p_i b_i+|\alpha|\ge s\right).$$
		To see this, take the class of $Q$ modulo $\FF_{Y,\mathbf p}^{s}$ and expand its base coefficients in $w$ modulo the indicated monomial ideals.
		In the fiberwise polynomial algebra, the surviving monomials $w^{\mathbf b}\eta^\alpha,$ with $\sum_i p_i b_i+|\alpha|<s,$ form a finite local basis over the ring of holomorphic functions of $z$ tensored with $\C[\xi].$
		The truncated class is zero at every stalk along $Z_{Y,k},$ so comparison of its polynomial coefficients in $\xi$ forces all these surviving coefficients to vanish.
		The resulting base monomial ideals are finitely generated, so the coefficients $q_{\alpha\beta}$ can be decomposed using finitely many generators without changing the jet exponents.
		Thus, $Q$ is a finite sum of terms $$\left(\prod_iw_i^{b_i}\right)A(w,z)\eta^\alpha\xi^\beta,\ \ \ \ \sum_i p_i b_i+|\alpha|\ge s,$$ each still of source weight $m,$ with $A$ holomorphic and depending only on the base coordinates.
		These coefficients remain holomorphic and bounded along the arc, so this finite expression can be evaluated even though the normal logarithmic jet coordinates have poles.
		The vanishing order of each such term is at least \begin{eqnarray*}\sum_i p_i'b_i-\sum_{i,h}h\beta_{i,h}&=&\sum_i p_i'b_i+\sum_{a,h}h\alpha_{a,h}-m\\&\ge&\sum_i p_i b_i+|\alpha|-m.\end{eqnarray*}
		The condition $\sum_i p_i b_i+|\alpha|\ge s$ bounds the last expression from below by $s-m.$
		Cancellation can only increase this order, proving \eqref{eq:order}.
		This includes an identically zero evaluation, whose order is $+\infty$ by convention.
		Taking the nonnegative part of $m-s$ and using the ceiling bound gives the stated pole estimate.
		The condition that the center lies in $Y^\circ$ is essential: residual logarithmic derivatives need not be holomorphic at a deeper stratum.
	\end{proof}
	~
	
	\begin{proof}[Proof of Theorem~\ref{thm:base}]~
		
		~
		
		Put $M=\mathscr O_{X_k(D)}(m)\otimes\pi_{0,k}^*H^{-a}.$
		Noetherianity gives positive integers $e_1,\ldots,e_v$ such that $$\bigcap_{\ell=1}^v\Bs|M^{e_\ell}|=\mathbf B(M).$$
		Taking a common multiple $q$ of their exponents yields $$\Bs|M^q|=\mathbf B(M)$$ as closed sets.
		Indeed, raising a section of $M^{e_\ell}$ to the power $q/e_\ell$ gives $\Bs|M^q|\subset\Bs|M^{e_\ell}|,$ while the reverse inclusion from $\mathbf B(M)$ follows from its definition.
		This exponent depends only on $M,$ not on the curve or its contacts.
		The direct-image identification \cite[(1.13)]{BD}, together with the projection formula, identifies sections of $M^q$ with invariant logarithmic jet differentials of weight $qm$ and twist $H^{-qa}.$
		The base locus hypothesis supplies a nonzero such differential $Q_0$ by choosing a regular jet over $X\setminus D$ and a section nonvanishing at its tower point.
		For every positive-dimensional stratum $Y,$ take $L=H^{qa},$ so that $\ell_Y=qa h_Y>0.$
		If $\mu_Y(\mathbf1)=0,$ Theorem~\ref{thm:mult} would force $Q_0=0,$ a contradiction.
		Hence, $\mu_Y(\mathbf1)>0.$
		By definition, this means that either $\rho_Y>0$ or $-\nu_{Y,i}>0$ for some $i,$ so $\mu_Y(\mathbf p)>0$ for every finite positive contact vector and all ratios in the statement are defined.
		Now fix a nonconstant entire curve $f$ not contained in $D.$
		Choose $z_0\in\C$ with $f(z_0)\notin D$ and $f'(z_0)\ne0,$ which is possible because the boundary preimage and the critical set are discrete.
		The point of $X_k(D)$ represented by $j_kf(z_0)$ is regular and lies outside the right side of \eqref{eq:base}, and hence outside $\Bs|M^q|.$
		Choosing a section of $M^q$ nonvanishing at this point gives, under the same direct-image identification, an invariant jet differential $$Q\in H^0\bigl(X,E_{k,qm}\Omega_X^1(\log D)\otimes H^{-qa}\bigr),\ \ \ \ Q(j_kf)\not\equiv0.$$
		In fact, this choice ensures $Q(j_kf(z_0))\ne0.$
		
		Fix $z_*\in f^{-1}(D),$ let $Y^\circ$ be the exact stratum containing $f(z_*),$ and write $J$ for its branch set and $p_i=\ord_{z_*}(w_i\circ f)$ for $i\in J.$
		Write $M_0=qm$ and define $$s_i=\left\lceil qa\theta_i(p_i)\right\rceil,\ \ \ \ s_Y=\left\lceil qa\theta_Y(\mathbf p)\right\rceil\ \ \ \ (\dim Y>0).$$
		Since $\sigma_{D_i}(H^{qa},p_i)\ge s_i$ and the filtration is decreasing, Theorem~\ref{thm:mult} gives membership in $\FF_{D_i,p_i}^{s_i}$ along the whole $Z_{D_i,k},$ including points over crossings.
		To impose these conditions simultaneously, use a common adapted coordinate chart and a local holomorphic frame of $H^{-qa}$ to write a finite expansion of source weight $M_0$ $$Q=\sum_{\alpha,\beta}A_{\alpha\beta}(w,z)\eta^\alpha\xi^\beta,$$ whose coefficients depend only on the base coordinates.
		For the component $D_i,$ the variables $\xi_{j,h}$ with $j\ne i$ and all the variables $\eta_{a,h}$ are tangential jet variables and have filtration degree one, while the variables $\xi_{i,h}$ have filtration degree zero.
		For a fixed full jet monomial, put $$j_i=\sum_{a,h}\alpha_{a,h}+\sum_{\substack{j\in J\\j\ne i}}\sum_h\beta_{j,h},\ \ \ \ d_i=\left\lceil\frac{(s_i-j_i)_+}{p_i}\right\rceil.$$
		The coefficient calculation in the proof of Theorem~\ref{thm:mult}, applied to the component filtration along the whole $Z_{D_i,k},$ shows that $w_i^{d_i}$ divides $A_{\alpha\beta}.$
		The functions $w_i$ are distinct prime parameters in the regular local coordinate ring, so these simultaneous divisibilities imply that $\prod_{i\in J}w_i^{d_i}$ divides the same coefficient $A_{\alpha\beta}.$
		Thus,  $Q$ has a common finite expansion into terms $$\left(\prod_{i\in J}w_i^{b_i}\right)A(w,z)\eta^\alpha\xi^\beta,\ \ \ \ p_i b_i+j_i\ge s_i\ \ (i\in J),$$ with $b_i=d_i,$ holomorphic base coefficients $A,$ and source weight $M_0$ in every term.
		For such a term, let $$M_i=\sum_hh\beta_{i,h}$$ be the source weight of the logarithmic derivatives belonging to branch $i.$
		Since every other jet factor has source weight at least one, $j_i\le M_0-M_i,$ and consequently $$M_i-p_i b_i\le M_0-s_i.$$
		Along the curve, the tangential jet variables for $Y^\circ$ are holomorphic and the normal logarithmic jet variable of order $h$ has pole order at most $h.$
		The signed pole order of the term is therefore at most $$\sum_{i\in J}(M_i-p_i b_i)\le |J|M_0-\sum_{i\in J}s_i,$$ and source homogeneity also bounds it by $M_0.$
		Cancellation can only increase the vanishing order, so both signed bounds also hold for the entire evaluation $Q(j_kf).$
		If $\dim Y>0,$ Theorem~\ref{thm:mult} gives the additional signed bound $M_0-s_Y$ for this evaluation at the exact stratum.
		In the following minimum, the last entry is omitted when $\dim Y=0.$
		Here, the signed pole order is $-\ord_{z_*}Q(j_kf),$ and the actual pole order is its nonnegative part. Taking the nonnegative part after the minimum, and using $M_0=qa c$ and the ceiling inequalities, gives $$\operatorname{pole}_{z_*}Q(j_kf)\le\left(\min\left\{M_0,\ |J|M_0-\sum_{i\in J}s_i,\ M_0-s_Y\right\}\right)_+\le qa\,b_Y(\mathbf p).$$ 
		Since $Q(j_kf)$ is holomorphic away from $f^{-1}(D),$ these pointwise estimates prove $$P_Q\le qa\,\mathcal B_f.$$
		Apply Lemma~\ref{lem:analytic} with $L=H^{qa}$ and its induced metric, and discard the nonnegative term $N(r,R_Q),$ to obtain $$qa\,T_f(r,H)\le N(r,P_Q)+S_f(r)\le qa\,N(r,\mathcal B_f)+S_f(r)\ \ \ \ \|.$$
		Dividing by $qa$ proves \eqref{eq:uniform}.
		For each curve the chosen section is fixed at all contacts, so varying contact orders introduce no new error constants or exceptional sets.
		The normalized coefficient involves $M_0/(qa)=m/a$ and the ratios in \eqref{eq:coefficient}, and is therefore independent of the chosen power $q.$
	\end{proof}
	~
	
	\begin{proof}[Proof of Theorem~\ref{thm:degeneracy}]~
		
		~
		
		Put $M=\mathscr O_{X_k(D)}(m)\otimes\pi_{0,k}^*H^{-a}$ and choose the exponent $q$ once so that $$\Bs|M^q|=\mathbf B(M),$$ as in the proof of Theorem~\ref{thm:base}.
		The open set $U=X\setminus(D\cup Z)$ is nonempty because $X$ is irreducible and $D,Z$ are proper closed subsets.
		Choosing a regular jet over $U$ and a section of $M^q$ nonvanishing at its tower point gives, by the direct-image identification, a nonzero global invariant jet differential $Q_0$ of weight $qm$ and twist $H^{-qa}.$
		For every positive-dimensional stratum $Y,$ including those contained in $Z,$ one has $\ell_Y=qa h_Y>0$ for $L=H^{qa}.$
		If $\mu_Y(\mathbf1)=0,$ Theorem~\ref{thm:mult} would force $Q_0=0,$ a contradiction.
		As in the proof of Theorem~\ref{thm:base}, it follows that $\mu_Y(\mathbf p)>0$ for every finite positive contact vector, so all required ratios are defined.
		Now let $f$ be a nonconstant entire curve not contained in $D\cup Z,$ and choose $z_0\in f^{-1}(U)$ with $f'(z_0)\ne0.$
		The set $f^{-1}(U)$ is nonempty and open, and the critical set of $f$ is discrete, so such a point exists.
		The corresponding regular tower point lies outside $\Bs|M^q|,$ so a section of this same power gives an invariant differential $$Q\in H^0\bigl(X,E_{k,qm}\Omega_X^1(\log D)\otimes H^{-qa}\bigr),\ \ \ \ Q(j_kf(z_0))\ne0.$$
		This choice of $Q$ is fixed along the whole curve.
		The local calculation in the proof of Theorem~\ref{thm:base} gives $$P_Q\le qa\mathcal B_f.$$
		The same local bounds apply at source points mapping into $D\cap Z,$ since they hold for the selected global differential $Q$ at every boundary contact.
		Applying Lemma~\ref{lem:analytic} and dividing by $qa$ proves the asserted second main theorem.
		
		For a positive-dimensional stratum, the definition \eqref{eq:mu} can be written as $$\mu_Y(\mathbf p)=\max\left\{0,\rho_Y,\frac{(-\nu_{Y,i})_+}{p_i}:i\in J\right\}.$$
		Each term depending on a contact variable is nonincreasing in that variable, so $\mu_Y$ is nonincreasing and the positive ratio $\theta_Y(\mathbf p)=h_Y/\mu_Y(\mathbf p)$ is nondecreasing.
		The same argument applies to the component ratios $\theta_i,$ and \eqref{eq:coefficient} shows that $b_Y$ is nonincreasing in every contact variable, also when $\dim Y=0$ and the $\theta_Y$ term is absent.
		Let $f$ now be a nonconstant entire curve not contained in $D$ with minimum contacts $\mathbf r,$ and suppose that $f(\C)\not\subset Z.$
		Then, $f^{-1}(D)$ and $f^{-1}(Z)$ are proper analytic subsets of $\C,$ hence discrete, and their union cannot cover $\C.$
		Therefore, $f(\C)\not\subset D\cup Z,$ and the first part supplies a differential $Q$ of weight $qm$ and twist $H^{-qa}$ with $Q(j_kf)\not\equiv0$ and $P_Q\le qa\mathcal B_f.$
		At a point $z_*\in f^{-1}(D),$ let $Y^\circ$ be the exact stratum containing $f(z_*)$ with branch set $J,$ and let $\mathbf p$ be the vector of actual contact orders.
		The minimum contact condition gives $p_i\ge r_i$ for $i\in J,$ so monotonicity and \eqref{eq:absorption} yield $$b_Y(\mathbf p)\le b_Y(\mathbf r_J)\le\sum_{i\in J}\lambda_i r_i\le\sum_{i\in J}\lambda_i p_i.$$
		The last sum is the coefficient of $f^*(\sum_i\lambda_iD_i)$ at $z_*,$ so these inequalities at all boundary contacts give $$P_Q\le qa\mathcal B_f\le f^*\left(qa\sum_i\lambda_iD_i\right).$$
		Set $L=H^{qa}$ and $R=qa\sum_i\lambda_iD_i.$
		Then $R$ is an effective rational divisor supported on $D,$ the chosen $Q$ is invariant of positive weight $qm,$ and its evaluation is not identically zero.
		The residual class is $$L-R=qa\left(H-\sum_i\lambda_iD_i\right),$$ which is ample by \eqref{eq:absorption}.
		All hypotheses of Lemma~\ref{lem:absorption} are therefore satisfied, and that lemma gives a contradiction.
		Thus, every such curve has image contained in the fixed set $Z.$
		This set is independent of the curve, and when $Z=\varnothing$ no such nonconstant entire curve exists.
	\end{proof}

	\begin{corollary}[Smooth boundary in arbitrary dimension]\label{cor:smooth}
		Suppose  that $D$ is smooth and irreducible and \eqref{eq:base} holds.
		Use the fixed triple $k,m,a$ from that hypothesis and write $c=m/a.$
		Choose an ample polarization $A$ on $D$ and put $$h=\deg_A(H|_D),\ \ \ \ \rho=\mu_{\max,A}(\Omega_D^1),\ \ \ \ \nu=\deg_A N_{D/X},\ \ \ \ \mu(p)=\max\{\rho,-\nu/p,0\}.$$
		Then,  $h>0$ and $\mu(p)>0$ for every integer $p\ge1.$
		For every nonconstant entire curve $f:\C\to X$ not contained in $D,$ the following inequalities hold
		\begin{equation}\label{eq:smooth}
			T_f(r,H)\le N\left(r,\sum_{z\in f^{-1}D}\left(c-\frac{h}{\mu(p_z)}\right)_+[z]\right)+S_f(r) \ \ \    \|,
		\end{equation}
		where $p_z=\ord_zf^*D.$
		In particular,  one has $$T_f(r,H)\le\left(c-\frac{h}{\mu(1)}\right)_+\Nbar_f(r,D)+S_f(r) \ \ \  \|.$$
		If all contacts have order at least an integer $p_0\ge1,$ one can replace $\mu(1)$ in the last formula by $\mu(p_0).$
	\end{corollary}
	
	\begin{proof}
		Since $n\ge2,$ the only closed boundary stratum is $Y=D,$ and it has positive dimension.
		Its residual boundary is empty, so $\Omega_D^1(\log B_D)=\Omega_D^1,$ and its normal bundle is $N_{D/X}=\OO_X(D)|_D.$
		Thus,  the numerical data of this stratum are precisely $h,\rho,\nu,\mu(p)$ as defined above.
		The ampleness of $H$ gives $h>0,$ and Theorem~\ref{thm:base} gives $\mu(p)>0$ for every integer $p\ge1.$
		There is only one branch, so the component and stratum contributions in \eqref{eq:coefficient} coincide, and $$b_D(p)=\left(\min\left\{c,c-\frac{h}{\mu(p)},c-\frac{h}{\mu(p)}\right\}\right)_+=\left(c-\frac{h}{\mu(p)}\right)_+.$$
		Theorem~\ref{thm:base} therefore gives \eqref{eq:smooth}.
		The identity $\mu(p)=\max\{\rho,0,(-\nu)_+/p\}$ shows that $\mu(p)$ is nonincreasing.
		Hence, for integers $p\ge p_0\ge1,$ we have $0<\mu(p)\le\mu(p_0)$ and $$\left(c-\frac{h}{\mu(p)}\right)_+\le\left(c-\frac{h}{\mu(p_0)}\right)_+.$$
		Finally, $(f^*D)_{\mathrm{red}}=\sum_{z\in f^{-1}D}[z].$ Thus,  for a curve whose contacts all have order at least $p_0,$ applying this coefficient inequality at each contact point and using the monotonicity of $N(r,\cdot),$ it  gives the asserted bound with $\mu(p_0),$ and in particular the bound with $\mu(1).$
	\end{proof}
	
	Formula~\eqref{eq:smooth} uses the actual contacts, while a common lower bound $p_0$ gives the uniform coefficient stated above.
	If $\nu\ge0,$ the function $\mu(p)$ is constant.
	If $\nu<0,$ then $\mu(p)=\max\{\max(\rho,0),(-\nu)/p\}.$
	When $\rho>0,$ this function decreases until it reaches $\rho;$ when $\rho\le0,$ it decreases to zero without attaining zero at a finite $p.$
	The coefficient decreases while $\mu(p)$ decreases and the coefficient remains positive.
	
	\begin{corollary}[A positive-coefficient degeneracy threshold]\label{cor:positive}
		Suppose that \eqref{eq:baseZ} holds for a fixed closed algebraic subset $Z\subsetneq X,$ possibly empty, and $D$ is smooth and irreducible.
		Here,  $c=m/a$ uses the triple in \eqref{eq:baseZ}, and $h,\rho,\nu,\mu$ are defined as in Corollary~\ref{cor:smooth}.
		Under this hypothesis, $h>0$ and $\mu(p)>0$ for every integer $p\ge1.$
		Let $d>0$ be rational and assume that $dH-D$ is nef.
		For an integer $p_0\ge1,$ put $\beta(p_0)=(c-h/\mu(p_0))_+.$
		If $$d\beta(p_0)<p_0,$$ then every nonconstant entire curve not contained in $D$ whose contacts have order at least $p_0$ has image contained in $Z.$
	\end{corollary}
	
	\begin{proof}
		Since $n\ge2,$ the only boundary stratum is $Y=D,$ and it has positive dimension and empty residual boundary.
		The ampleness of $H$ gives $h>0,$ and the first assertion of Theorem~\ref{thm:degeneracy} gives $\mu(p)>0$ for every integer $p\ge1$ under \eqref{eq:baseZ}.
		The single-branch form of \eqref{eq:coefficient} is $$b_D(p_0)=\left(c-\frac{h}{\mu(p_0)}\right)_+=\beta(p_0).$$
		The strict numerical hypothesis gives $\beta(p_0)/p_0<1/d,$ so we can choose a rational number $\lambda$ with $$0\le\frac{\beta(p_0)}{p_0}\le\lambda<\frac1d.$$
		Consequently, $b_D(p_0)\le\lambda p_0,$ which is the only stratum inequality required in \eqref{eq:absorption} for the minimum contact order $p_0.$
		Moreover, $$H-\lambda D=(1-\lambda d)H+\lambda(dH-D).$$
		Here, $1-\lambda d>0$ and $\lambda\ge0,$ so the first summand is an ample rational class and the second is a nef rational class.
		Their sum is therefore ample, and Theorem~\ref{thm:degeneracy} gives the asserted conclusion.
	\end{proof}
	
	The numerical equivalence $D\equiv dH$ is a special case of the nef domination used here.
	Unlike that special case, the nef-domination hypothesis also allows negative normal degree, so the contact dependence of $\mu(p_0)$ can contribute to the criterion.
	When $\nu\ge0,$ the threshold is $p_0>d(c-h/\rho)_+$ rather than the full-weight threshold $p_0>dc$ for the same base locus data.
	This comparison does not assert optimality over other constructions of jet differentials.
	
	\begin{corollary}[A conditional orbifold Green-Griffiths-Lang criterion]\label{cor:ggl}
		Fix integers $r_i\ge1$ and set $$\Delta_{\mathbf r}=\sum_{i=1}^t\left(1-\frac1{r_i}\right)D_i,$$ using the minimum-contact convention for orbifold entire curves \cite[Definition 1.3]{CDDR}.
		Suppose that $K_X+\Delta_{\mathbf r}$ is big and \eqref{eq:baseZ} holds for a fixed closed algebraic subset $Z\subsetneq X,$ possibly empty.
		Suppose  that \eqref{eq:absorption} holds for this contact vector $\mathbf r$ and the coefficients determined by the same triple $k,m,a.$
		Then, every nonconstant orbifold entire curve of $(X,\Delta_{\mathbf r})$ has image in the fixed proper algebraic subset $D\cup Z.$
		More precisely, every such curve not contained in $D$ has image in  $Z$ from \eqref{eq:baseZ}.
	\end{corollary}
	
	\begin{proof}
		Let $f:\C\to X$ be a nonconstant orbifold entire curve of $(X,\Delta_{\mathbf r}).$
		By definition, $f(\C)\not\subset\operatorname{Supp}\Delta_{\mathbf r}.$
		If $f(\C)\subset D,$ its image already lies in $D\cup Z.$
		Suppose now that $f(\C)\not\subset D,$ so that $f^*D_i$ is a well-defined effective divisor for every $i.$
		For $r_i>1,$ the orbifold condition gives $\ord_zf^*D_i\ge r_i$ at every $z\in f^{-1}(D_i).$
		For $r_i=1,$ the component $D_i$ does not occur in $\operatorname{Supp}\Delta_{\mathbf r},$ but the same inequality holds automatically because the order at every actual intersection is a positive integer.
		Thus, $f$ has minimum contacts $\mathbf r$ along the full divisor $D,$ and Theorem~\ref{thm:degeneracy} gives $f(\C)\subset Z.$
		Finally, $D$ and $Z$ are fixed proper closed algebraic subsets of the irreducible variety $X,$ so $D\cup Z\subsetneq X$ is a proper closed algebraic subset independent of $f.$
	\end{proof}
	
	The general-type condition in this corollary places the conclusion in the setting of the orbifold Green-Griffiths-Lang conjecture \cite[Conjecture 1.5]{CDDR}; it is not used to deduce the jet or numerical hypotheses.
	In particular, the conclusion is uniform algebraic degeneracy, stronger than the assertion that each individual curve has a proper Zariski closure.
	For the logarithmic complement, omission of $D$ already makes every evaluation holomorphic, and \eqref{eq:baseZ} gives the established uniform degeneracy consequence of jet vanishing \cite{Dem,BD}.
	The additional content here is a sufficient numerical test for finite minimum contacts with positive corrected pole coefficients.
	For a fixed effective boundary class $R=\sum_i\lambda_iD_i$ with $H-R$ ample, the inequalities in \eqref{eq:absorption} specify which finite orbifold multiplicities are allowed.
	This does not prove the conjecture from bigness alone, nor establish a new unconditional family without separately verifying the base locus condition.
	
	\begin{corollary}[A numerical ramification obstruction]\label{cor:ramification}
		Let $f:\C\to X$ be a nonconstant entire curve not contained in $D.$
		Assume that \eqref{eq:base} is satisfied.
		Use the fixed triple $k,m,a$ in this hypothesis and put $c=m/a.$
		Then, we have $\mathcal B_f\ne0.$
		In particular, suppose that $D$ is smooth and irreducible, and use the notation of Corollary~\ref{cor:smooth}.
		If $\nu<0$ and an integer $p_0\ge1$ satisfies $$c\rho\le h,\ \ \ \ p_0h\ge c(-\nu),$$ then no such curve can have all its contacts with $D$ of order at least $p_0.$
	\end{corollary}
	
	\begin{proof}
		If $\mathcal B_f=0,$ the construction in the proof of Theorem~\ref{thm:base} supplies a positive integer $q$ and a global invariant differential $Q$ of weight $M=qm$ and twist $H^{-qa}$ such that $$Q(j_kf)\not\equiv0,\ \ \ \ P_Q\le qa\,\mathcal B_f=0.$$
		Since $P_Q$ is effective, it follows that $P_Q=0.$
		With $L=H^{qa}$ and $R=0,$ the residual class $L-R$ is ample, so Lemma~\ref{lem:absorption} gives a contradiction.
		This proves $\mathcal B_f\ne0.$
		For the smooth-boundary assertion, suppose that all contacts have order at least $p_0.$
		For every integer $p\ge p_0,$ the two numerical conditions give $$\rho\le\frac{h}{c},\ \ \ \ 0<\frac{-\nu}{p}\le\frac{-\nu}{p_0}\le\frac{h}{c}.$$
		Consequently, $$0<\mu(p)=\max\left\{\rho,\frac{-\nu}{p},0\right\}\le\frac{h}{c}.$$
		By the single-branch form of \eqref{eq:coefficient}, the coefficient at any actual contact of order $p\ge p_0$ is $$b_D(p)=\left(c-\frac{h}{\mu(p)}\right)_+=0.$$
		Thus, we have  $\mathcal B_f=0,$ also when the contact set is empty, contradicting with the first part.
	\end{proof}
	
	Under the numerical hypotheses above, for a smooth irreducible boundary the corollary excludes nonconstant orbifold entire curves of $(X,(1-1/p_0)D)$ whose images are not contained in $D,$ with contacts interpreted as minimum multiplicities. For $p_0\ge2,$ the restriction on the image is part of the orbifold definition; for $p_0=1,$ it must be retained explicitly.
	The curves omitting $D$ give the established consequence of the stable base locus hypothesis; the numerical test also treats curves meeting $D.$
	
	There is also a useful consistency condition when $D$ is smooth and irreducible with  nonnegative normal degree.
	For a line bundle $L$ over $X,$ one puts $\ell=\deg_A(L|_D).$
	If $\nu\ge0,$ the finite source weight bounds the tangential exponent $j$ in Lemma~\ref{lem:filtration} by $m.$
	Thus, every section $Q\in H^0(X,\GG_{k,m}\Omega_X^1(\log D)\otimes L^{-1})$ with $\ell>m\max\{\rho,0\}$ is zero: all graded factors have negative maximal slope at every filtration degree.
	Equivalently, the constrained feasible set defining $\sigma_D^{[m]}$ is empty.
	Consequently, \eqref{eq:base} forces $\rho>0$ and $c\rho\ge h$ in this case.
	For a hypersurface boundary this gives a necessary numerical condition on admissible weight-to-twist ratios; it does not supply the required sections.
	
	For comparison with optimization over the same construction, let $\mathcal A$ be a nonempty collection of triples satisfying \eqref{eq:base} for a fixed pair with smooth irreducible boundary and fixed $H.$
	Set $c_{\mathcal A}=\inf_{\mathcal A}m/a$ and $\theta=h/\mu(1).$
	For each $\varepsilon>0,$ choosing a single triple sufficiently close to this infimum gives $$T_f(r,H)\le(c_{\mathcal A}-\theta+\varepsilon)_+\Nbar_f(r,D)+S_f(r) \ \ \      \|.$$
	The subtraction therefore persists when the same collection of weight-to-twist ratios is optimized.
	For each $\varepsilon$ the proof uses a single choice of numerical data and a single exceptional set.
	The ratio $h_Y/\mu_Y(\mathbf p)$ can likewise be optimized over the polarization when it varies with the ample class.
	Rescaling a polarization multiplies $h_Y,\rho_Y,$ and every $\nu_{Y,i}$ by the same factor, so only its ray matters for these ratios.
	The ratio is independent of the polarization on a curve or on a stratum of Picard number one.
	
	\subsection{Proofs of Theorems~\ref{prop:fiber} and \ref{prop:cross} and examples}~
	
	\begin{proof}[Proof of Theorem~\ref{prop:fiber}]~
		
		~
		
		We work with the algebraic jet space and ideal sheaves as in Lemma~\ref{lem:filtration}.
		The regular embedding gives $$\II_Y^r/\II_Y^{r+1}\simeq\Sym^r(\II_Y/\II_Y^2).$$
		For an integer $r\ge0$ satisfying $r\eta_{\xi_0}<\ell_Y,$ Lemma~\ref{lem:slope} gives $$\mu_{\max,A_Y}\bigl(\Sym^rE_{\xi_0}\otimes L^{-1}|_Y\bigr)\le r\eta_{\xi_0}-\ell_Y<0,$$ where $r=0$ uses $\Sym^0E_{\xi_0}=\OO_Y.$
		Hence $$H^0\bigl(Y,\Sym^rE_{\xi_0}\otimes L^{-1}|_Y\bigr)=0.$$
		The product projection $Z_{Y,k}=Y\times\mathbb A^{|J|k}\to\mathbb A^{|J|k}$ is projective, of finite presentation, and flat, and the vector bundle $$\Sym^r(\II_Y/\II_Y^2)\otimes\pi_Y^*(L^{-1}|_Y)$$ is flat over the parameter space.
		Upper semicontinuity \cite[Tag 0BDN]{Stacks} therefore gives an open neighborhood $U_r$ of $\xi_0$ on which $$H^0\bigl(Y,\Sym^rE_\xi\otimes L^{-1}|_Y\bigr)=0\ \ \ \ (\xi\in U_r).$$
		The open set $U_r$ may depend on $r.$
		Suppose inductively that $Q\in\II_Y^r\otimes\pi^*L^{-1}$ along $Z_{Y,k},$ starting with $r=0.$
		Its graded image is a global section of the vector bundle above, and its restriction to every fiber over $U_r$ is zero.
		It therefore vanishes on the dense open subset $Y\times U_r$ of the integral variety $Z_{Y,k},$ and local freeness forces it to vanish everywhere.
		The exact sequence of consecutive ideal powers then gives $Q\in\II_Y^{r+1}\otimes\pi^*L^{-1}$ along $Z_{Y,k}.$
		If $\eta_{\xi_0}>0,$ iteration for $$0\le r<\left\lceil\frac{\ell_Y}{\eta_{\xi_0}}\right\rceil$$ proves the asserted ideal-power membership.
		If $\eta_{\xi_0}=0,$ the same argument applies to every $r\ge0.$
		At every $\zeta\in Z_{Y,k},$ Krull intersection gives $$\bigcap_{r\ge0}\bigl(\II_Y^r\otimes\pi^*L^{-1}\bigr)_\zeta=0,$$ since $\II_{Y,\zeta}$ is contained in the maximal ideal of the Noetherian local ring $\mathscr O_{\JJ_k,\zeta}.$
		Thus, $Q$ vanishes on a nonempty open subset of the irreducible jet space $\JJ_k,$ and hence $Q=0.$
		
		For the partition, consider, for each branch $i\in J,$ the polynomial map $$(\xi_{i,1},\ldots,\xi_{i,k})\longmapsto\bigl(B_1(\xi_i),\ldots,B_k(\xi_i)\bigr).$$
		This is a triangular polynomial automorphism with diagonal coefficients one, since $B_h(\xi_i)-\xi_{i,h}$ depends only on $\xi_{i,1},\ldots,\xi_{i,h-1}.$
		We can therefore choose $\xi_0,$ so that $$\left.B_h(\xi_i)\right|_{\xi=\xi_0}=\begin{cases}1,&i\in G_h,\\0,&i\notin G_h.\end{cases}$$
		Write $[w_i]$ and $[\eta_{a,h}]$ for the induced local frame elements of $E_{\xi_0}.$
		In the notation of the proof of Lemma~\ref{lem:filtration}, the transition laws on this fiber become \begin{eqnarray*}{[w_i']}&=&u_i|_Y\,[w_i],\\ {[\eta'_{a,h}]}&=&\sum_bA_{ab}[\eta_{b,h}]+\sum_{i\in G_h}C_{ai}[w_i].\end{eqnarray*}
		Thus, the transition matrix separates into blocks indexed by $h,$ with the $h$th block generated by $[\eta_{a,h}]$ and $[w_i]$ for $i\in G_h.$
		On $X_{G_h},$ a local frame of $\Omega_{X_{G_h}}^1(\log R_{G_h})$ consists of $dw_i$ for $i\in G_h$ and the forms $\vartheta_a$ used in that lemma.
		After restriction of this bundle to $Y,$ its transition matrix is exactly the $h$th block above under $[w_i]\leftrightarrow dw_i$ and $[\eta_{a,h}]\leftrightarrow\vartheta_a.$
		This also covers empty groups and proves $$E_{\xi_0}\simeq\bigoplus_{h=1}^kF_{G_h}.$$
		The simple normal crossing assumption gives the logarithmic conormal exact sequence $$0\longrightarrow\bigoplus_{i\in G}N_i^*\longrightarrow F_G\longrightarrow\Omega_Y^1(\log B_Y)\longrightarrow0,$$ with $F_{\varnothing}=\Omega_Y^1(\log B_Y)$ when $G=\varnothing.$
		Since $\deg_{A_Y}N_i^*=-\nu_{Y,i},$ Lemma~\ref{lem:slope} gives $$\mu_{\max,A_Y}(F_G)\le\max\bigl(\{\rho_Y\}\cup\{-\nu_{Y,i}:i\in G\}\bigr)\le\mu_Y(\mathbf1).$$
		Taking the maximum over the direct summands and with zero yields $$\eta_{\xi_0}=\max\bigl\{0,\mu_{\max,A_Y}(F_{G_h}):1\le h\le k\bigr\}=\eta_{\mathcal P}\le\mu_Y(\mathbf1).$$
	\end{proof}
	
	For an arc centered in $Y^\circ$ with contacts at least $\mathbf p,$ the two independent multiplicity estimates give $$\ord_0Q(j_k\gamma)\ge\max\left\{\sigma_Y(L,\mathbf p),\left\lceil\frac{\ell_Y}{\eta_{\mathcal P}}\right\rceil\right\}-m,\ \ \ \ \eta_{\mathcal P}>0.$$
	Thus, the thresholds are combined by taking their maximum.
	In \eqref{eq:coefficient} this replaces $\theta_Y(\mathbf p)$ by $$\max\left\{\theta_Y(\mathbf p),\frac{h_Y}{\eta_{\mathcal P}}\right\}.$$
	The estimate from the full conormal bundle remains unweighted; the contact dependence comes from Theorem~\ref{thm:mult}.
	For one smooth irreducible boundary and $k\ge2,$ choose $B_1(\xi)=1$ and $B_h(\xi)=0$ for $h\ge2.$
	The resulting conormal fiber is $\Omega_X^1|_D\oplus(\Omega_D^1)^{\oplus(k-1)}.$
	Using the notation of Corollary~\ref{cor:smooth}, the conormal exact sequence $$0\longrightarrow N_{D/X}^*\longrightarrow\Omega_X^1|_D\longrightarrow\Omega_D^1\longrightarrow0$$ and Lemma~\ref{lem:slope} give $$-\nu\le\mu_{\max,A}(\Omega_X^1|_D)\le\max\{-\nu,\rho\}.$$
	Since $k\ge2,$ the direct sum above contains $\Omega_D^1$ as a direct summand, so $$\eta_\xi=\max\{0,\mu_{\max,A}(\Omega_X^1|_D),\rho\}=\max\{0,-\nu,\rho\}=\mu_D(1).$$
	Thus, this fiber supplies no additional unweighted improvement.
	In order one, taking $B_1(\xi)=1$ makes the conormal fiber $\Omega_X^1|_D,$ recovering a direct application of the general conormal principle \cite{Pet}.
	At $\xi=0,$ the conormal fiber is $N_{D/X}^*\oplus\Omega_D^1.$
	
	~
	
	\begin{proof}[Proof of Theorem~\ref{prop:cross}]~
		
		~
		
		Choose adapted coordinates $x,y$ with $D_i=(x=0)$ and $D_j=(y=0),$ and write $\xi_{i,h}=(\log x)^{(h)}$ and $\xi_{j,h}=(\log y)^{(h)}.$
		Then $$\II_i=(x,\xi_{j,1},\ldots,\xi_{j,k}),\ \ \ \ \II_j=(y,\xi_{i,1},\ldots,\xi_{i,k}).$$
		Write $$Q=\sum_{\alpha,\beta}a_{\alpha\beta}(x,y)\prod_{h=1}^k\xi_{i,h}^{\alpha_h}\xi_{j,h}^{\beta_h}.$$
		The sum is finite, its coefficients are holomorphic, and each term satisfies $$\sum_{h=1}^kh(\alpha_h+\beta_h)=m.$$
		For a fixed term, one  puts $A=\sum_h\alpha_h$ and $B=\sum_h\beta_h,$ so $A+B\le m.$
		The two ideal-power conditions give $$x^{(s_i-B)_+}\mid a_{\alpha\beta},\ \ \ \ y^{(s_j-A)_+}\mid a_{\alpha\beta}.$$
		Since $x,y$ are distinct prime parameters, the product of these powers divides the same coefficient.
		Factoring out this product expresses $Q$ as a finite sum of terms $x^by^da(x,y)\prod_h\xi_{i,h}^{\alpha_h}\xi_{j,h}^{\beta_h}$ with $a$ holomorphic, $b,d\in\mathbb Z_{\ge0},$ and $$A+B\le m,\ \ \ \ b+B\ge s_i,\ \ \ \ d+A\ge s_j.$$
		For fixed $b,d,$ these inequalities admit nonnegative integers $A,B$ if and only if $$(s_i-b)_++(s_j-d)_+\le m.$$
		Using the identity $$u_++v_+=\max\{0,u,v,u+v\},$$ this is equivalent to $$b\ge b_0,\ \ \ \ d\ge d_0,\ \ \ \ b+d\ge\delta.$$
		The definitions also give $$\delta-b_0-d_0=\min\big\{s_i,s_j,m,(s_i+s_j-m)_+\big\}\ge0.$$
		For every such pair $b,d,$ we have \begin{eqnarray*}p_ib+p_jd&=&p_ib_0+p_jd_0+p_i(b-b_0)+p_j(d-d_0)\\ &\ge&p_ib_0+p_jd_0+\min(p_i,p_j)(b+d-b_0-d_0)\\ &\ge&C.\end{eqnarray*}
		Equality is attained at the feasible integer pair $$(b_*,d_*)=\begin{cases}(\delta-d_0,d_0),&p_i\le p_j,\\(b_0,\delta-b_0),&p_i>p_j.\end{cases}$$
		Along $\gamma,$ each logarithmic jet variable of order $h$ has pole order at most $h,$ while $x^by^d$ has vanishing order $p_ib+p_jd$ and the remaining coefficient is holomorphic.
		Thus, every term of $Q(j_k\gamma)$ has vanishing order at least $p_ib+p_jd-m\ge C-m.$
		The same bound holds for their finite sum, including an identically zero evaluation, so its pole order is at most $(m-C)_+.$
		To prove optimality for the stated local information, take the minimizing pair $(b_*,d_*)$ above and put $$B=(s_i-b_*)_+,\ \ \ \ A=m-B.$$
		The feasibility condition gives $0\le B\le m$ and $A\ge(s_j-d_*)_+.$
		Consequently, the local differential $$Q_*=x^{b_*}y^{d_*}\xi_{i,1}^{A}\xi_{j,1}^{B}$$ belongs to $\II_i^{s_i}\cap\II_j^{s_j}$ and has source weight exactly $m$ for every $k\ge1.$
		For the arc $\gamma_*(t)=(t^{p_i},t^{p_j})$ in these coordinates, we have $$Q_*(j_k\gamma_*)=p_i^{A}p_j^{B}t^{C-m}.$$
		The leading constant is nonzero, so its pole order is exactly $(m-C)_+,$ also when $C\ge m$ and this order is zero.
		This proves the asserted local optimality.
	\end{proof}
	
	On a surface, the slope data become intersection numbers.
	For $D_i$ a boundary component, $$\rho_i=(K_X+D)\cdot D_i,\ \ \ \ \nu_i=D_i^2,\ \ \ \ h_i=H\cdot D_i,\ \ \ \ \mu_i(p)=\max\{\rho_i,-\nu_i/p,0\}.$$
	At a smooth boundary point the stable coefficient is $(c-h_i/\mu_i(p))_+.$
	At a double point the component part of \eqref{eq:coefficient} is $$\left(\min\left\{c,2c-\frac{h_i}{\mu_i(p_i)}-\frac{h_j}{\mu_j(p_j)}\right\}\right)_+.$$
	For unit contacts and finite-section thresholds $s_i,s_j,$ the signed calculation gives $(\min\{m,2m-s_i-s_j\})_+.$
	This retains the signs before taking a positive part, and can be stronger than $\min\{m,(m-s_i)_++(m-s_j)_+\}.$
	
	Theorem~\ref{prop:cross} supplies the additional contact improvement used below.
	
	If $s_i,s_j\le m,$ Theorem~\ref{prop:cross} reduces to $$\operatorname{pole}_0Q(j_k\gamma)\le\bigl(m-\min(p_i,p_j)(s_i+s_j-m)_+\bigr)_+.$$
	In the stable-base-locus setting take unweighted $\theta_i=h_i/\mu_i(1)$ and define $$b_0'=(\theta_i-c)_+,\ \ \ \ d_0'=(\theta_j-c)_+,\ \ \ \ \delta'=(\theta_i+\theta_j-c)_+,$$ $$C'=p_i b_0'+p_jd_0'+\min(p_i,p_j)(\delta'-b_0'-d_0').$$
	The double-point coefficient can additionally be replaced by its minimum with $(c-C')_+.$
	Apply Theorem~\ref{prop:cross} with source weight $qm$ and the unweighted integer thresholds $s_\ell=\lceil qa\theta_\ell\rceil$ for $\ell\in\{i,j\}.$
	After division by $qa,$ the continuous relaxation retains the thresholds $s_i/(qa)$ and $s_j/(qa).$
	Lowering these thresholds to $\theta_i$ and $\theta_j$ enlarges the feasible set and gives the minimum $C'.$
	Thus the integer minimum in Theorem~\ref{prop:cross} is at least $qaC'.$
	
	\begin{example}[A strict contact gain at a double point]
		Take $k\ge1,$ $m=5,$ $s_i=s_j=3,$ and contacts $(p_i,p_j)=(2,3).$
		Theorem~\ref{prop:cross} gives $b_0=d_0=0,$ $\delta=1,$ and $C=2,$ hence a pole bound of $3.$
		If one uses only the lower bounds $p_i,p_j\ge1,$ the same two ideal conditions give the larger pole bound $$\big(m-(s_i+s_j-m)_+\big)_+=4.$$
		In local coordinates $D_i=(x=0)$ and $D_j=(y=0),$ the two ideals are $$\II_i=(x,\xi_{j,1},\ldots,\xi_{j,k}),\ \ \ \ \II_j=(y,\xi_{i,1},\ldots,\xi_{i,k}).$$
		Consider the monomial $$Q=x\,\xi_{i,1}^{3}\xi_{j,1}^{2}.$$
		Since $x\xi_{j,1}^{2}\in\II_i^3$ and $\xi_{i,1}^{3}\in\II_j^3,$ we have $Q\in\II_i^3\cap\II_j^3.$
		The factor $x$ has source weight zero, so $Q$ has source weight $3+2=5.$
		Along the arc $\gamma(t)=(t^2,t^3),$ one has $$\xi_{i,1}(j_k\gamma)=\frac2t,\ \ \ \ \xi_{j,1}(j_k\gamma)=\frac3t,$$ and hence $$Q(j_k\gamma)=t^2\left(\frac2t\right)^3\left(\frac3t\right)^2=72t^{-3}.$$
		Its pole order is therefore exactly $3,$ so the contact-dependent bound is attained.
		The example proves local sharpness under the stated ideal constraints; extending the monomial to a global section is a separate problem.
	\end{example}
	
	For a smooth irreducible boundary curve $D$ on a surface, let $L$ be a line bundle with $\ell=L\cdot D>0,$ and put $\kappa=2g(D)-2$ and $\nu=D^2.$
	The integer threshold \eqref{eq:sigma} is $$\sigma_p=\min\big\{j+pb: \   j,b\in\mathbb Z_{\ge0},\ j\kappa-b\nu\ge\ell\big\}.$$
	When $\nu<0$ and $\kappa>0,$ it has the finite expression $$\sigma_p=\min_{0\le j\le\lceil\ell/\kappa\rceil}\left[j+p\left\lceil\frac{(\ell-j\kappa)_+}{-\nu}\right\rceil\right].$$
	When $\nu<0$ and $\kappa\le0,$ it is $$\sigma_p=p\left\lceil\frac{\ell}{-\nu}\right\rceil.$$
	This can improve the simpler ceiling $\lceil p\ell/(-\nu)\rceil.$
	When $\kappa\le0$ and $\nu<0,$ this is already a direct consequence of the standard logarithmic jet filtration \cite[Proposition 2.4(c)]{CDDR} and the restricted residue sequence \cite[Fact 11.4]{GKKP}.
	Indeed, these give $\mu_{\max}(\GG_{k,m}\Omega_X^1(\log D)|_D)\le0.$
	For every integer $b<\ell/(-\nu),$ restriction therefore gives $$H^0\bigl(D,\GG_{k,m}\Omega_X^1(\log D)|_D\otimes L^{-1}|_D\otimes N_{D/X}^{-b}\bigr)=0.$$
	With $s=\lceil\ell/(-\nu)\rceil,$ successive divisor restriction sequences give $$Q\in H^0\bigl(X,\GG_{k,m}\Omega_X^1(\log D)\otimes L^{-1}\otimes\mathscr O_X(-sD)\bigr).$$
	Pulling back this fixed boundary factor gives the displayed contact saving.
	
	\begin{example}[Negative normal degree with nonzero sections]\label{ex:ruled}
		Let $C$ be a smooth projective curve of genus two and let $M$ be a line bundle of degree seven on $C.$
		Compactify the total space of $M$ to a ruled surface $\pi:X\to C,$ and let $D_0$ be its infinity section.
		Then,  $N_{D_0/X}\simeq M^{-1}$ and $D_0^2=-7.$
		Choose three general sections of $M,$ with graphs $S_1,S_2,S_3,$ and set $$D=D_0+S_1+S_2+S_3.$$
		The bundle $M$ separates the  $1$-jets due to  that  $\deg M=7>2g(C),$ thus the graphs may be chosen with transverse pairwise intersections and no triple intersections.
		Each graph satisfies $\mathscr O_X(S_i)\simeq\mathscr O_X(D_0)\otimes\pi^*M.$
		They are disjoint from $D_0,$ and their pairwise intersection numbers are seven.
		Thus, $D$ is SNC and its residual boundary on $D_0$ is empty.
		
		Choose a line bundle $B$ of degree seventeen on $C$ and put $$H=\mathscr O_X(D_0)\otimes\pi^*B.$$
		One has $$H^2=27,\ \ \ \ H\cdot D_0=10,$$ and $H$ has degree one on each fiber.
		For every other horizontal irreducible curve $T,$ one has $$H\cdot T=D_0\cdot T+17\deg(\pi|_T)>0,$$ since $D_0\cdot T\ge0.$
		Thus, $H$ is ample by the Nakai--Moishezon criterion \cite[Chapter V, Theorem 1.10]{Hart}.
		Along $D_0,$ the data are $$\rho=2,\ \ \ \ \nu=-7,\ \ \ \ h=H\cdot D_0=10,\ \ \ \ \mu(p)=\max\{2,7/p\}.$$
		For twist $H^{-1}$ and an integer $p\ge1,$ the integer threshold is $$\sigma_p=\min\big\{j+pb: \    j,b\in\mathbb Z_{\ge0},\ 2j+7b\ge10\big\}=\min\{5,2p,p+2\}.$$
		For $b=0,1,2,$ the least possible $j$ is respectively $5,2,0,$ while $b\ge3$ gives cost at least $3p>2p.$
		For source weight $m\ge5,$ all three minimizing alternatives are allowed and $\sigma_{D_0}^{[m]}(H,p)=\sigma_p.$
		The thresholds compare as follows.
		\begin{center}
			\begin{tabular}{lcccc}
				\toprule
				Contact order $p$ \ &\  $1$ \ & \ $2$ \ & \ $3$ \ & $p\ge4$ \\
				\midrule
				Integer threshold $\sigma_p$ \ & \ $2$ \ & \ $4$ \ & $5$ & $5$ \\
				Slope ceiling $\lceil10/\mu(p)\rceil$ \ & \ $2$ \ &\  $3$ \ & $5$ & $5$ \\
				\bottomrule
			\end{tabular}
		\end{center}
		At a double contact the integer threshold therefore gives an additional unit of pole saving beyond the slope ceiling.
		The saturation for $p\ge3$ reflects the tangential alternative $j=5,$ which has cost independent of contact order.
		
		We next verify that this pair admits nonzero negatively twisted sections.
		Write $F$ for a fiber and $E=\Omega_X^1(\log D).$
		Here, $D_0\cdot F=1,$ $F^2=0$ and the boundary has numerical class $D\equiv4D_0+21F.$
		The standard ruled-surface formulas and the logarithmic Chern formula give $$K_X\equiv-2D_0-5F,\ \ \ \ K_X+D\equiv2D_0+16F,\ \ \ \ c_1(E)^2=36,\ \ \ \ c_2(E)=25.$$
		Here, $e$ denotes the topological Euler characteristic.
		For the last equality, we have $e(X)=-4,$ every boundary component has Euler number $-2,$ and inclusion-exclusion at the twenty-one transverse double points gives $$c_2(E)=e(X\setminus D)=-4-4(-2)+21=25.$$
		Riemann--Roch \cite[Appendix A, Section 4]{Hart} consequently gives $$\chi\bigl(X,\Sym^mE\otimes H^{-1}\bigr)=\frac{c_1(E)^2-c_2(E)}{6}m^3+O(m^2)=\frac{11}{6}m^3+O(m^2).$$
		Let $U\subset C$ be the complement of the projections of the twenty-one crossings.
		For every fiber $F$ over $U,$ the four boundary points are distinct, so the relative logarithmic cotangent sequence restricts to $$0\longrightarrow\mathscr O_F\longrightarrow E|_F\longrightarrow\mathscr O_F(2)\longrightarrow0.$$
		The quotient is $\Omega_F^1(\log(D|_F))\simeq\mathscr O_F(2),$ because $F\simeq\PP^1$ and $D|_F$ consists of four reduced points.
		Dualizing this sequence and tensoring by $\mathscr O_F(-1)$ gives an extension of $\mathscr O_F(-1)$ by $\mathscr O_F(-3),$ and therefore $$H^0\bigl(F,(E|_F)^*\otimes\mathscr O_F(-1)\bigr)=0.$$
		The splitting theorem on $F\simeq\PP^1$ then implies that every summand of $E|_F$ has nonnegative degree.
		Since $(H\otimes K_X)|_F\simeq\mathscr O_F(-1),$ every summand of $\bigl(\Sym^mE^*\otimes H\otimes K_X\bigr)|_F$ has negative degree, so this bundle has no nonzero section on $F.$
		Any global section consequently vanishes on the dense open subset $\pi^{-1}(U)$ and is zero.
		Serre duality therefore gives $$H^2\bigl(X,\Sym^mE\otimes H^{-1}\bigr)=H^0\bigl(X,\Sym^mE^*\otimes H\otimes K_X\bigr)^*=0.$$
		Consequently $$h^0\bigl(X,\Sym^mE\otimes H^{-1}\bigr)\ge\chi\bigl(X,\Sym^mE\otimes H^{-1}\bigr)>0$$ for all sufficiently large $m.$
		The canonical inclusion $$\Sym^mE\otimes H^{-1}\hookrightarrow\GG_{k,m}\Omega_X^1(\log D)\otimes H^{-1}\ \ \ \ (k\ge1)$$ is obtained from homogeneous polynomials involving only first-order jet variables.
		Thus, the integer and contact estimates apply to nonzero global sections on this pair, independently of a stable base locus condition.
		The example compares the two universal vanishing bounds on a nonzero section space.
		Whether a global section attains either bound is a separate question.
		In fact, this pair fails \eqref{eq:base} at order one.
		Write $s_i\in H^0(C,M)$ for the section with graph $S_i$ and put $\varphi=(s_3-s_1)/(s_2-s_1).$
		The absence of triple intersections gives $\deg\varphi=7,$ and pairwise transversality makes $\varphi:C\to\PP^1$ unramified over $0,1,\infty.$
		By Riemann--Hurwitz, the ramification divisor of $\varphi$ has degree sixteen, so there is a critical point $b\in C$ with $\varphi(b)\notin\{0,1,\infty\}.$
		Near $b,$ the fiber coordinate $v=(u-s_1)/(s_2-s_1),$ where $u$ is the original fiber coordinate in a local trivialization, puts the four boundary sections at $0,1,\varphi,\infty.$
		Along the fiber $F_b,$ differentiation in a base coordinate $t$ at fixed $v$ gives a nowhere-vanishing logarithmic tangent vector.
		Indeed, near the moving section one may hold $v-\varphi(t)$ fixed instead, and these two horizontal vectors agree on $F_b$ because $\varphi'(b)=0.$
		Evaluating any section of $\Sym^m E\otimes H^{-a}$ on this vector gives a global section of $H^{-a}|_{F_b},$ which vanishes since $\deg(H^{-a}|_{F_b})=-a<0.$
		The corresponding tangent directions over $F_b\setminus D$ therefore belong to the base locus of every positive power in \eqref{eq:base} at $k=1.$
	\end{example}
	
	For a smooth plane curve of degree $d\ge4$ and $H=\mathscr O(1),$ adjunction gives $$\kappa=d(d-3),\ \ \ \ \nu=d^2,\ \ \ \ h=d.$$
	Thus, every section of weight $m$ and twist $\mathscr O(-a)$ has pole order at most $$\left(m-\left\lceil\frac{a}{d-3}\right\rceil\right)_+,$$ and Corollary~\ref{cor:smooth} gives the coefficient $(c-1/(d-3))_+$ under \eqref{eq:base}.
	For an SNC plane boundary of total degree $d\ge4,$ put $\theta=1/(d-3);$ every component then has $\theta_i=\theta.$
	The nonnegative normal degrees imply $c\ge\theta$ under \eqref{eq:base}, by the source-weight obstruction above.
	Let $\Sigma$ be the finite crossing set and define $\Nbar_f(r,\Sigma)=N(r,\sum_{f(z)\in\Sigma}[z]).$
	The general stratum bound becomes $$T_f(r,\OO(1))\le(c-\theta)\Nbar_f(r,D)+\min\{\theta,c-\theta\}\Nbar_f(r,\Sigma)+S_f(r) \ \ \    \|.$$
	Indeed, the coefficient is $c-\theta$ at a smooth boundary point and $\min\{c,2c-2\theta\}$ at a double point.
	The crossing correction is therefore $\theta$ for $c\ge2\theta,$ and $c-\theta$ for $\theta\le c\le2\theta.$
	Here, $\Nbar_f(r,D)$ counts a crossing once, in contrast to the sum of the reduced counting functions of the components.
	
	\begin{corollary}[Plane SNC boundaries and unequal orbifold multiplicities]\label{cor:planeorb}
		Let $D=\sum_iD_i\subset\PP^2$ be a simple-normal-crossing divisor with smooth components of degrees $d_i$ and total degree $d\ge4.$
		Assume \eqref{eq:baseZ} with $H=\mathscr O(1)$ and a fixed closed algebraic subset $Z\subsetneq\PP^2,$ possibly empty, and put $$\theta=\frac1{d-3},\ \ \ \ \beta=c-\theta.$$
		Here, $c=m/a$ uses the fixed triple in \eqref{eq:baseZ}.
		Then $\beta\ge0,$ and every nonconstant entire curve not contained in $D\cup Z$ satisfies $$T_f(r,\mathscr O(1))\le\beta\sum_i\Nbar_f(r,D_i)+S_f(r) \ \ \     \|.$$
		If integers $r_i\ge1$ satisfy \begin{equation}\label{eq:planeorb}\beta\sum_i\frac{d_i}{r_i}<1,\end{equation} all nonconstant entire curves with these minimum contacts and not contained in $D$ have images in the same set $Z.$
		Set $$\Delta_{\mathbf r}=\sum_i\left(1-\frac1{r_i}\right)D_i.$$
		If also $\sum_i d_i/r_i<d-3,$ then $(\PP^2,\Delta_{\mathbf r})$ is of general type and every nonconstant orbifold entire curve, understood in the minimum-contact sense, has image contained in the fixed proper algebraic subset $D\cup Z.$
	\end{corollary}
	
	\begin{proof}
		For each component $D_i,$ one has $h_i=d_i$ and $\nu_i=d_i^2.$
		The residual boundary $(D-D_i)|_{D_i}$ has degree $d_i(d-d_i),$ so adjunction gives $$\rho_i=d_i(d_i-3)+d_i(d-d_i)=d_i(d-3)>0.$$
		Consequently, $\mu_i(p_i)=\rho_i$ and $\theta_i(p_i)=\theta$ for every integer $p_i\ge1.$
		As in the proof of Theorem~\ref{thm:degeneracy}, hypothesis \eqref{eq:baseZ} supplies a positive integer $q$ and a nonzero global invariant differential $Q_0$ of source weight $qm$ and twist $H^{-qa}.$
		If $c<\theta,$ then for every component $D_i$ and all integers $0\le j\le qm$ and $b\ge0,$ one has $$j\rho_i-b\nu_i\le qm\,d_i(d-3)<qa\,d_i=\deg(H^{qa}|_{D_i}).$$
		The constrained feasible set in the finite source weight refinement of Theorem~\ref{thm:mult} is therefore empty, forcing $Q_0=0,$ a contradiction.
		Thus, $c\ge\theta$ and $\beta\ge0.$
		At a smooth point of $D_i$ the coefficient in \eqref{eq:coefficient} is $\beta,$ and at a double point it is $\min\{c,2\beta\}\le2\beta,$ independently of the contact orders.
		A plane SNC divisor has no triple points, so for every curve in the first assertion we obtain the effective-divisor inequality $$\mathcal B_f\le\beta\sum_i(f^*D_i)_{\mathrm{red}}.$$
		At a source point mapping to a crossing, the sum on the right has coefficient $2\beta,$ since that point is counted once for each of the two components.
		Applying the second main inequality in Theorem~\ref{thm:degeneracy} and taking counting functions proves the asserted estimate.
		For the uniform conclusion, put $\lambda_i=\beta/r_i.$
		These are nonnegative rational numbers, and the single-component and double-point stratum conditions are respectively 
		$$\beta=\lambda_i r_i,\ \ \ \  \min\{c,2\beta\}\le2\beta=\lambda_i r_i+\lambda_j r_j. $$
		They cover every stratum inequality in \eqref{eq:absorption}, while $$H-\sum_i\lambda_iD_i\equiv\left(1-\beta\sum_i\frac{d_i}{r_i}\right)H$$ is ample by \eqref{eq:planeorb}.
		Theorem~\ref{thm:degeneracy} therefore gives the stated conclusion with the same fixed set $Z.$
		Finally, $$K_{\PP^2}+\Delta_{\mathbf r}\equiv\left(d-3-\sum_i\frac{d_i}{r_i}\right)H.$$
		Under the additional inequality this class is ample, hence big, and Corollary~\ref{cor:ggl} gives the asserted orbifold conclusion.
	\end{proof}
	
	For a common minimum contact $p,$ condition \eqref{eq:planeorb} is $$p>d\left(c-\frac1{d-3}\right),$$ compared with $p>dc$ from the full-weight estimate with the same data.
	The strict inequality is essential for the ample residual class; no conclusion at equality is claimed.
	Corollary~\ref{cor:planeorb} allows unequal multiplicities and treats crossings without assuming that the curve avoids the crossing set.
	
	\begin{example}[Hypersurface slopes and the necessity of higher jets]
		Let $X=\PP^n$ with $n\ge2,$ let $D$ be a smooth hypersurface of degree $d\ge n+2,$ and take $H=\mathscr O_{\PP^n}(1).$
		Use $A=\mathscr O_D(1)$ as polarization.
		The canonical bundle $K_D=\mathscr O_D(d-n-1)$ is ample.
		The negative K\"ahler--Einstein metric \cite{Yau} has K\"ahler class proportional to $c_1(K_D)$ and makes the cotangent bundle
		 Hermitian-Einstein, hence polystable \cite[Theorem 5.8.3]{Kob}.
		Since $c_1(K_D)=(d-n-1)c_1(A),$ this polystability holds with respect to the chosen polarization $A.$
		For $n=2$ the same slope computation is simply the degree of a line bundle.
		In particular, $\mu_{\max,A}(\Omega_D^1)=\mu_A(\Omega_D^1).$
		Using $A^{n-1}=d$ and $N_{D/\PP^n}\simeq\mathscr O_D(d),$ we obtain $$h=d,\ \ \ \ \nu=d^2,\ \ \ \ \rho=\frac{d(d-n-1)}{n-1},\ \ \ \ \theta=\frac{n-1}{d-n-1}.$$
		Since $\rho,\nu>0,$ we have $\mu(p)=\rho$ and $h/\mu(p)=\theta$ for every positive integer $p.$
		Let $$Q\in H^0\bigl(\PP^n,\GG_{k,m}\Omega_{\PP^n}^1(\log D)\otimes\OO(-a)\bigr),\ \ \ \ k,m,a\ge1.$$
		For every holomorphic arc $\gamma$ centered at a point of $D$ and not contained in $D,$ Theorem~\ref{thm:mult} gives $$\operatorname{pole}_0Q(j_k\gamma)\le\left(m-\left\lceil\frac{a(n-1)}{d-n-1}\right\rceil\right)_+.$$
		Whenever \eqref{eq:base} holds, with $c=m/a$ for its fixed triple, the resulting counting coefficient is $$\left(c-\frac{n-1}{d-n-1}\right)_+.$$
		The correction is therefore explicit once the numerical data $k,m,a$ of a section construction are specified.
		More generally, under \eqref{eq:baseZ}, use $c=m/a$ from its fixed triple.
		Since $dH-D\equiv0$ is nef, Corollary~\ref{cor:positive} confines every nonconstant entire curve not contained in $D$ and with all contacts of order at least a fixed integer $p\ge1$ to $Z$ whenever $$p>d\left(c-\frac{n-1}{d-n-1}\right)_+.$$
		The associated orbifold canonical class is $$K_{\PP^n}+(1-1/p)D\equiv\left(d-n-1-\frac d p\right)H,$$ so the additional inequality $p>d/(d-n-1)$ places this conclusion in the general-type setting of Corollary~\ref{cor:ggl}.
		
		There is also a classical first-order obstruction in this setting.
		For $k=1,$ choose $\xi$ with $B_1(\xi)=1$ in Theorem~\ref{prop:fiber}, so that $E_\xi\simeq\Omega_{\PP^n}^1|_D.$
		The restricted Euler sequence is $$0\longrightarrow\Omega_{\PP^n}^1|_D\longrightarrow\OO_D(-1)^{\oplus(n+1)}\longrightarrow\OO_D\longrightarrow0.$$
		Consequently, $\mu_{\max,A}(E_\xi)\le-d<0$ and $\eta_\xi=0.$
		For $L=\OO(a),$ one has $\ell_D=ad>0,$ and Theorem~\ref{prop:fiber} therefore gives $$H^0\bigl(\PP^n,\Sym^m\Omega_{\PP^n}^1(\log D)\otimes\OO(-a)\bigr)=0,\ \ \ \ m,a\ge1.$$
		Thus, nonzero sections with this twist require $k\ge2.$
		For $k\ge2,$ choose a normal parameter with $B_1(\xi)=1$ and $B_h(\xi)=0$ for $h\ge2.$
		Theorem~\ref{prop:fiber} then gives the full decomposition $$E_\xi\simeq\Omega_{\PP^n}^1|_D\oplus(\Omega_D^1)^{\oplus(k-1)}.$$
		Its maximal slope is $\rho>0,$ and hence $\eta_\xi=\rho.$
		Higher jet order therefore removes the obstruction exhibited by this conormal fiber, although existence of the sections remains a global problem; see the constructions in \cite{BD,HHMX}.
	\end{example}
	
	The next three examples verify the stable base locus condition \eqref{eq:base} for concrete logarithmic surface pairs.
	They provide instances of the geometric hypothesis in Theorem~\ref{thm:base}, and hence of \eqref{eq:baseZ} with $Z=\varnothing$ in Theorem~\ref{thm:degeneracy}.
	In each case, a fixed positive power of the relevant line bundle has no base points over the complement of the boundary.
	
	\begin{example}[A del Pezzo surface with ten boundary curves]\label{ex:base-delpezzo}
		Let $X$ be the blowup of $\PP^2$ at four points $p_1,\ldots,p_4$ with no three collinear, and let $D$ be the sum of the four exceptional curves and the strict transforms of the six lines joining pairs of these points.
		Then $X$ is a del Pezzo surface of degree five, $H=-K_X$ is very ample and $D\sim2H,  H^2=5.$
		The divisor $D$ consists of ten smooth rational curves with simple normal crossings.
		We prove \eqref{eq:base} with $(k,m,a)=(1,6,1).$
		
		Write $\tau:X\to\PP^2$ for the blowdown, $\Lambda=\tau^*\OO_{\PP^2}(1),$ and $E=\sum_{i=1}^4E_i$ for its exceptional divisor.
		Number the six lines as $L_j=(\ell_j=0),$ denote their strict transforms by $\widetilde L_j,$ and set $I_i=\{j:p_i\in L_j\}.$
		Each $I_i$ has three elements, and $$\widetilde L_j\sim\Lambda-\sum_{i:p_i\in L_j}E_i,\ \ \ \ \sum_{j=1}^6\widetilde L_j=\tau^*\left(\sum_{j=1}^6L_j\right)-3E.$$
		Since $H\sim3\Lambda-E,$ this gives $$D=\sum_j\widetilde L_j+E\sim6\Lambda-2E\sim2H.$$
		
		The classical Hirzebruch--Kummer construction of exponent five gives a finite Galois cover $\pi:Y\to X$ with group $G\simeq(\mathbb Z/5\mathbb Z)^5,$ branched with index five along all ten components of $D.$
		The smooth projective surface $Y$ is a compact ball quotient; see \cite[Introduction]{BCcover}.
		Here, one takes the normalization of $X$ in the extension $$\C(\PP^2)\bigl((\ell_j/\ell_1)^{1/5}:2\le j\le6\bigr).$$
		A product of these five ratios can be a fifth power only when all its exponents are divisible by five, as follows by taking orders along $L_2,\ldots,L_6.$
		Thus, the extension has the stated Galois group.
		In additive notation, identify $$G=(\mathbb Z/5\mathbb Z)^6/\langle(1,\ldots,1)\rangle,$$ and let $g_j$ be the image of the $j$th standard basis vector, so that $\sum_jg_j=0.$
		For $j\in I_i$ and $j_0\notin I_i,$ the chosen extension contains $(\ell_j/\ell_{j_0})^{1/5},$ while $\ord_{E_i}(\ell_j/\ell_{j_0})=1.$
		The inertia along $\widetilde L_j$ is generated by $g_j,$ and the inertia along $E_i$ is generated by $\sum_{j\in I_i}g_j,$ which is nonzero and has order five.
		At a boundary crossing, the two inertia generators are either $g_j,g_k$ for $j\ne k,$ or $g_j,\sum_{k\in I_i}g_k$ for $j\in I_i,$ and are linearly independent over $\mathbb F_5.$
		For the second pair, a linear relation in this quotient lifts to a multiple of $(1,\ldots,1);$ examining first the coordinates outside $I_i,$ then those in $I_i\setminus\{j\},$ and finally the $j$th coordinate forces both coefficients to vanish.
		
		The ramification along $E_i$ can also be seen directly in local coordinates.
		Choose a reference line not passing through $p_i,$ which does not change the field extension, and coordinates in which the three incident lines have equations $u=0,$ $v=0,$ and $u-v=0.$
		After taking fifth roots of local holomorphic units, the local extension is obtained by adjoining fifth roots of $u,$ $v,$ and $u-v.$
		In the blowup chart $v=r,\ u=rt,$ it is equivalently obtained by adjoining fifth roots of $r,$ $t,$ and $t-1.$
		After these local choices, each connected local component of the normalization is described by $$r=z^5,\ \ \ \ t=\alpha^5,\ \ \ \ \beta^5=\alpha^5-1.$$
		The three original fifth roots are $z\alpha,z,z\beta.$
		Let $\zeta$ be a primitive fifth root of unity.
		Multiplying these roots by $\zeta^A,\zeta^B,\zeta^C,$ respectively, with $A,B,C\in\mathbb Z/5\mathbb Z,$ induces $$ (z,\alpha,\beta)\longmapsto(\zeta^Bz,\zeta^{A-B}\alpha,\zeta^{C-B}\beta).$$
		Such an automorphism acts trivially on the residue field $\C(\alpha,\beta)$ at the generic point of $(z=0)$ precisely when $A=B=C.$
		The inertia subgroup is therefore generated by simultaneous multiplication of all three original roots, which is the element $\sum_{j\in I_i}g_j.$
		The reduced preimage of $E_i=(r=0)$ is given by $z=0,$ and the normal parameter satisfies $r=z^5,$ so the ramification index is five.
		The induced map on the ramification curve is generically \'etale in the tangential direction, while the surface map remains ramified in the normal direction.
		At a point $(0,\alpha_0,\beta_0)$ above a smooth boundary point of $E_i,$ put $t_0=\alpha_0^5\notin\{0,1,\infty\}.$
		Since $\alpha_0\beta_0\ne0,$ the implicit function theorem shows that $\eta=\alpha^5-t_0$ is a local parameter on the ramification curve.
		With source coordinates $(z,\eta)$ and target coordinates $(r,t-t_0),$ the map becomes $(z,\eta)\mapsto(z^5,\eta).$
		
		Put $R=(\pi^{-1}D)_{\mathrm{red}}.$
		After suitable local analytic changes of coordinates, the models over a smooth boundary point and a crossing are respectively $(u,v)\mapsto(u^5,v)$ and $(u,v)\mapsto(u^5,v^5).$
		The reduced preimage $R$ is locally $(u=0)$ in the first case and $(uv=0)$ in the second.
		Consequently, $R$ is SNC and $$\pi^*\Omega_X^1(\log D)\simeq\Omega_Y^1(\log R).$$
		Since every component has ramification index five, the ramification formula gives $$K_Y\sim\pi^*K_X+4R,\ \ \ \ \pi^*D=5R.$$
		Consequently, $$K_Y^{\otimes5}\simeq\pi^*\bigl(K_X^{\otimes5}\otimes\OO_X(4D)\bigr)\simeq(\pi^*H)^{\otimes3}.$$
		
		Let $p:\PP(T_Y)\to Y$ parameterize tangent lines, so that $p_*\OO_{\PP(T_Y)}(r)=\Sym^r\Omega_Y^1$ for $r\ge0.$
		The ball metric makes $\OO_{\PP(T_Y)}(1)$ positive, and the line bundle $$\mathscr N=\OO_{\PP(T_Y)}(3)\otimes p^*K_Y^{-1}$$ is nef.
		For completeness, normalize the curvature tensor of the ball metric at a unitary frame as $$R_{i\bar j k\bar l}=-b(\delta_{ij}\delta_{kl}+\delta_{il}\delta_{kj}),\ \ \ \ b>0.$$
		At $(y,[v])\in\PP(T_Y),$ the induced curvature forms on a Chern-horizontal vector $u$ satisfy \begin{eqnarray*}\Theta_{\OO(1)}(u,\bar u)&=&b\left(\|u\|^2+\frac{|\langle u,v\rangle|^2}{\|v\|^2}\right),\\ \Theta_{K_Y}(u,\bar u)&=&3b\|u\|^2.\end{eqnarray*}
		Thus, the horizontal curvature of $\mathscr N$ is $3b|\langle u,v\rangle|^2/\|v\|^2\ge0,$ its vertical curvature is $3\omega_{\mathrm{FS}}>0,$ and the mixed terms vanish in this horizontal splitting.
		These are the standard projective-bundle curvature formulas; see \cite[proof of Lemma 1.3 and the discussion following Lemma 1.4]{BKT}.
		Set $$\mathscr A=\OO_{\PP(T_Y)}(6)\otimes p^*\pi^*H^{-1}.$$
		The preceding line-bundle identity yields $$\mathscr A^{\otimes3}\simeq\OO_{\PP(T_Y)}(3)\otimes\mathscr N^{\otimes5}.$$
		The right side is ample, so $\mathscr A$ is ample.
		
		Choose one integer $q\ge1$ such that $\mathscr A^q$ is globally generated, and put $N=|G|.$
		Fix $x\in X\setminus D$ and $0\ne v\in T_{X,x}.$
		Since $\pi$ is \'etale over $X\setminus D,$ the vector $v$ has one nonzero lift at each of the $N$ points of $\pi^{-1}(x).$
		Global generation and the fact that finitely many proper linear subspaces cannot cover a complex vector space allow us to choose $$P\in H^0\bigl(Y,\Sym^{6q}\Omega_Y^1\otimes\pi^*H^{-q}\bigr)$$ whose evaluations on all these lifted directions are nonzero.
		Using the natural $G$-equivariant inclusion $$\Omega_Y^1\hookrightarrow\Omega_Y^1(\log R)\simeq\pi^*\Omega_X^1(\log D),$$ the product $\prod_{g\in G}g^*P$ is a $G$-invariant section of $$\pi^*\bigl(\Sym^{6qN}\Omega_X^1(\log D)\otimes H^{-qN}\bigr).$$
		The projection formula and $(\pi_*\OO_Y)^G=\OO_X$ descend this product to a global section $Q$ on $X$ with $Q_x(v)\ne0.$
		The integers $q$ and $N$ have been fixed independently of $x$ and $v.$
		Since $X_1(D)=\PP(T_X(-\log D))$ and first-order invariant jet differentials are symmetric differentials, we obtain $$\Bs\left|\left(\OO_{X_1(D)}(6)\otimes\pi_{0,1}^*H^{-1}\right)^{qN}\right|\subset\pi_{0,1}^{-1}(D).$$
		This proves \eqref{eq:base} for the stated fixed triple, with $c=6.$
	\end{example}
	
	\begin{example}[Blowing up a boundary crossing]\label{ex:base-node}
		Let $(X,D,H)$ be the pair in Example~\ref{ex:base-delpezzo}, and let $\sigma:X'\to X$ be the blowup of a crossing point of $D,$ with exceptional curve $F.$
		Define $$D'=(\sigma^{-1}D)_{\mathrm{red}}=\sigma^*D-F,\ \ \ \ H'=2\sigma^*H-F.$$
		The divisor $D'$ is SNC, and $F$ meets the strict transforms of the two boundary branches at two distinct points.
		Since $H$ is very ample, the hyperplane sections through the center of the blowup generate $\sigma^*H-F,$ so this line bundle is nef.
		For every irreducible curve $C\ne F$  $$H'\cdot C=\sigma^*H\cdot C+(\sigma^*H-F)\cdot C>0.$$
		Moreover, $H'\cdot F=1$ and $(H')^2=19,$ so the Nakai--Moishezon criterion gives that $H'$ is ample.
		
		The natural logarithmic differential map is an isomorphism $$\sigma^*\Omega_X^1(\log D)\simeq\Omega_{X'}^1(\log D').$$
		Indeed, if $D=(xy=0)$ near the center, then the chart $x=u,\ y=uv$ has $D'=(uv=0)$ and $$\sigma^*(d\log x)=d\log u,\ \ \ \ \sigma^*(d\log y)=d\log u+d\log v.$$
		The other blowup chart gives the same invertibility.
		By Example~\ref{ex:base-delpezzo}, there is a fixed integer $s\ge1$ such that every nonzero direction over $X\setminus D$ has nonzero evaluation for some $$Q\in H^0\bigl(X,\Sym^{6s}\Omega_X^1(\log D)\otimes H^{-s}\bigr).$$
		Let $s_F$ be the canonical section of $\OO_{X'}(F).$
		Then $$Q'=s_F^s(\sigma^*Q)^2\in H^0\bigl(X',\Sym^{12s}\Omega_{X'}^1(\log D')\otimes(H')^{-s}\bigr),$$ because $(H')^{-s}=(\sigma^*H)^{-2s}\otimes\OO_{X'}(sF).$
		Over $X'\setminus D',$ the map $\sigma$ is an isomorphism and $s_F$ is nonzero, so these sections detect every nonzero tangent direction there.
		Consequently, $$\Bs\left|\left(\OO_{X'_1(D')}(12)\otimes\pi_{0,1}^*(H')^{-1}\right)^s\right|\subset\pi_{0,1}^{-1}(D'),$$ proving \eqref{eq:base} with $(k,m,a)=(1,12,1)$ and $c=12.$
	\end{example}
	
	\begin{example}[Blowing up a smooth boundary point]\label{ex:base-smooth}
		Let $(X,D,H)$ again be as in Example~\ref{ex:base-delpezzo}, choose a point $q\in D_i$ lying on no other boundary component, and let $\sigma:X''\to X$ be its blowup, with exceptional curve $F.$
		Put $$D''=(\sigma^{-1}D)_{\mathrm{red}}=\sigma^*D,\ \ \ \ H''=2\sigma^*H-F.$$
		Here, $F$ is included in the boundary and meets the strict transform of $D_i$ transversely at one point.
		So, $D''$ is SNC and the intersection argument in Example~\ref{ex:base-node} proves that $H''$ is ample.
		
		There is a natural regular map $$\sigma^*\Omega_X^1(\log D)\longrightarrow\Omega_{X''}^1(\log D''),$$ which is an isomorphism over $X''\setminus D''.$
		To check regularity along $F,$ choose local coordinates with $D=(x=0).$
		In the chart $x=u,\ y=uv,$ the new boundary is $(u=0)$ and $$\sigma^*(d\log x)=d\log u,\ \ \ \ \sigma^*(dy)=uv\,d\log u+u\,dv.$$
		In the chart $x=zt,\ y=t,$ the new boundary is $(zt=0)$ and $$\sigma^*(d\log x)=d\log z+d\log t,\ \ \ \ \sigma^*(dy)=t\,d\log t.$$
		The two matrices have determinants $u$ and $t,$ respectively, so the map is regular and degenerates along $F.$
		
		Take the same fixed exponent $s$ and the sections $Q$ supplied by Example~\ref{ex:base-delpezzo}, and view their pullbacks as logarithmic symmetric differentials through this natural map.
		With $s_F$ the canonical section of $\OO_{X''}(F),$ the sections $$Q''=s_F^s(\sigma^*Q)^2\in H^0\bigl(X'',\Sym^{12s}\Omega_{X''}^1(\log D'')\otimes(H'')^{-s}\bigr)$$ detect every nonzero tangent direction over $X''\setminus D''.$
		It follows that $$\Bs\left|\left(\OO_{X''_1(D'')}(12)\otimes\pi_{0,1}^*(H'')^{-1}\right)^s\right|\subset\pi_{0,1}^{-1}(D''),$$ so \eqref{eq:base} again holds with $(k,m,a)=(1,12,1)$ and $c=12.$
		The exceptional curve belongs to the boundary, and no additional algebraic exceptional set is needed.
	\end{example}
	
	Examples~\ref{ex:base-delpezzo}--\ref{ex:base-smooth} realize the stable base locus hypothesis  at 1-jets.
	The two blowups give different smooth SNC compactifications of the complement in Example~\ref{ex:base-delpezzo}.
	
	The following examples realize \eqref{eq:base} at higher orders and establish the minimum jet order in each case.
	We first provide  a common argument for general hypersurfaces with hyperplane sections as boundary.
	
	Let $n\ge2,$ let $X\subset\PP^{n+1}$ be a general smooth hypersurface of degree $d\ge2(n+1)^{2n+6},$ and set $H=\OO_X(1).$
	Let $B$ be the reduced union of $s$ hyperplanes in linear general position, with $1\le s\le n+2,$ chosen so that $X$ meets every stratum of $B$ transversely.
	Then $D=B|_X$ is reduced and SNC.
	By Deng \cite[Theorem C and Corollary D]{DengDT}, there is a fixed weight $m_0\ge1$ such that for each regular $n$-jet on $X$ some global section of $E_{n,m_0}\Omega_X^1\otimes H^{-1}$ has nonzero evaluation.
	The natural inclusion into logarithmic jet differentials preserves evaluation over $X\setminus D.$
	Indeed, an ordinary derivative normal to a boundary branch is a regular polynomial in the logarithmic jet coordinates by the identities $x^{(h)}=xB_h((\log x)',\ldots,(\log x)^{(h)}).$
	The inclusion respects changes of source parameter, so it preserves the invariant subsheaves.
	It follows that \eqref{eq:base} holds with $(k,m,a)=(n,m_0,1)$ and that the required inclusion already holds for the base locus of the line bundle itself.
	
	To exclude all lower orders, write $E=\Omega_X^1(\log D)$ and $V=\Omega_{\PP^{n+1}}^1(\log B)|_X.$
	The hyperplanes can be taken to be coordinate hyperplanes after a projective change of coordinates.
	In the logarithmic Euler sequence the surjection $\OO^{\oplus s}\oplus\OO(-1)^{\oplus(n+2-s)}\to\OO$ is given by contraction with the Euler vector field and splits through the first trivial summand.
	Thus $$V\simeq\OO_X^{\oplus(s-1)}\oplus\OO_X(-1)^{\oplus(n+2-s)}.$$
	The transversality to all boundary strata gives the exact logarithmic conormal sequence $$0\longrightarrow\OO_X(-d)\longrightarrow V\longrightarrow E\longrightarrow0.$$
	Locally the defining equation of $X$ can be completed to an adapted coordinate system together with the local equations of the boundary branches, so its differential generates the kernel, also at the crossings.
	For every $\ell\ge1,$ taking the degree-$\ell$ part of the symmetric algebra yields the exact sequence $$0\longrightarrow\OO_X(-d)\otimes\Sym^{\ell-1}V\longrightarrow\Sym^\ell V\longrightarrow\Sym^\ell E\longrightarrow0.$$
	For $\ell=0,$ one has $\Sym^0 E=\OO_X.$
	
	The logarithmic Green-Griffiths filtration \cite[Proposition 2.4(c)]{CDDR} at order $j$ has graded factors $\bigotimes_{h=1}^{j}\Sym^{\ell_h}E$ with $\sum_{h=1}^{j}h\ell_h=w,$ where $w$ is the weight.
	Tensoring the locally split resolutions just obtained and twisting by $H^{-a}$ for $a\ge1$ gives a resolution of each twisted factor of length at most $j$ by direct sums of negative line bundles $\OO_X(-b),$ with $b>0.$
	The hypersurface exact sequence in projective space gives $H^i(X,\OO_X(-b))=0$ for $0\le i<n$ and $b>0.$
	When $j<n,$ the long cohomology sequences of the resolution therefore give vanishing of $H^0$ for each twisted factor.
	The filtration and the inclusion of invariant jet differentials into Green-Griffiths jet differentials give $$H^0\bigl(X,E_{j,w}\Omega_X^1(\log D)\otimes H^{-a}\bigr)=0,\ \ \ \ 1\le j<n,\ \ \ \ w,a\ge1.$$
	This includes the weights and twists of every positive power in the definition of the stable base locus, so the minimum order for \eqref{eq:base} with this polarization is exactly $n.$
	
	\begin{example}[A surface with smooth boundary requiring 2-jets]\label{ex:base-order-two}
		Let $X\subset\PP^3$ be a general smooth surface of degree $d\ge2\cdot3^{10}=118098,$ let $D$ be a smooth hyperplane section, and take $H=\OO_X(1).$
		The argument above with $n=2$ and $s=1$ proves \eqref{eq:base} at order two and excludes it at order one.
	\end{example}
	
	\begin{example}[A surface with two boundary components requiring 2-jets]
		Let $X\subset\PP^3$ be a general smooth surface of degree $d\ge118098,$ and let $D=D_1+D_2$ be the sum of two general hyperplane sections.
		The components meet transversely at the $d$ points where their common line meets $X.$
		The argument above with $n=s=2$ shows that the minimum order for \eqref{eq:base} with $H=\OO_X(1)$ is two.
	\end{example}
	
	\begin{example}[A threefold requiring 3-jets]
		Let $X\subset\PP^4$ be a general smooth hypersurface of degree $d\ge2\cdot4^{12}=33554432,$ let $D$ be a smooth hyperplane section, and take $H=\OO_X(1).$
		The argument above with $n=3$ and $s=1$ proves \eqref{eq:base} at order three and gives vanishing for all negative twists by powers of $H$ at orders one and two.
	\end{example}
	
	\begin{example}[2-jets after blowing up a smooth boundary point]
		Take $(X,D,H)$ from Example~\ref{ex:base-order-two} and let $\sigma:X'\to X$ be the blowup of a point of $D,$ with exceptional curve $F.$
		Write $\widetilde D$ for the strict transform and set $$D'=\sigma^*D=\widetilde D+F,\ \ \ \ H'=2\sigma^*H-F.$$
		The divisor $D'$ is SNC.
		The bundle $\sigma^*H-F$ is globally generated by hyperplanes through the center because $H$ is very ample.
		Thus, $H'$ has positive degree on every curve other than $F$ by its decomposition as $\sigma^*H+(\sigma^*H-F).$
		Since $H'\cdot F=1$ and $(H')^2=4d-1>0,$ the Nakai--Moishezon criterion gives ampleness.
		
		Let $m_0$ be the fixed weight supplied by the construction on $X,$ and let $s_F$ be the canonical section of $\OO_{X'}(F).$
		For each regular two-jet over $X'\setminus D',$ choose a section $Q\in H^0(X,E_{2,m_0}\Omega_X^1\otimes H^{-1})$ that detects its image under $\sigma.$
		Pullback of this ordinary jet differential by the holomorphic map $\sigma,$ followed by the inclusion into logarithmic jet differentials, gives $$s_F(\sigma^*Q)^2\in H^0\bigl(X',E_{2,2m_0}\Omega_{X'}^1(\log D')\otimes(H')^{-1}\bigr).$$
		Its evaluation is nonzero on the chosen jet because $\sigma$ is an isomorphism over the complement of the boundary and $s_F$ does not vanish there.
		This proves \eqref{eq:base} with $(k,m,a)=(2,2m_0,1)$ already for the first power.
		
		More generally, any global section $Q'$ of $E_{k,w}\Omega_{X'}^1(\log D')\otimes(H')^{-a}$ for $k,w,a\ge1$ can be written as $Q'=s_F^a\sigma^*Q$ with $$Q\in H^0\big(X,E_{k,w}\Omega_X^1(\log D)\otimes H^{-2a}\big).$$
		On $X'\setminus F,$ dividing by $s_F^a$ gives a section on the complement of the blowup center in $X,$ and Hartogs extension for locally free sheaves extends it over that point.
		The logarithmic pullback is regular by the chain rule in the blowup charts used below, so the identity extends across $F.$
		In particular, all such sections vanish over $F,$ and order-one sections vanish identically by the preceding lower-order calculation on $X.$
		
		At an exact boundary point on $F\setminus\widetilde D,$ one has $h_F=1,$ $\rho_F=-1,$ and $\nu_F=-1,$ since $F\simeq\PP^1$ and meets the other boundary component at one point.
		Thus, $\mu_F(p)=1/p,$ so \eqref{eq:coefficient} gives $b_F(p)=(2m_0-p)_+$ for a contact of order $p.$
		The factorization through $X$ gives a stronger local pole bound, also at the crossing.
		On $D$ the slope is $d(d-3)$ by adjunction and the normal degree is $d,$ while the degree of $H^{2a}|_D$ is $2ad.$
		Consequently, Theorem~\ref{thm:mult} with unit contact weight gives membership of $Q$ in the $\tau$th power of the pure normal jet ideal, with $\tau=\lceil2a/(d-3)\rceil.$
		In local coordinates $(x,y)$ with $D=(x=0),$ this ideal is generated by $x,y',\ldots,y^{(k)}.$
		In the chart $x=u,\ y=uv,$ the logarithmic pullback of each of these generators is divisible by $u,$ since every ordinary derivative of $u$ is $u$ times a polynomial in its logarithmic jet coordinates.
		In the chart $x=zt,\ y=t,$ the same generators pull back to multiples of $t.$
		The logarithmic derivatives of $x$ are those of $u$ in the first chart and the sums of those of $z$ and $t$ in the second, so the pullback is regular in both charts.
		The charts cover $F,$ with $u=0$ and $t=0$ its respective equations, and the second chart contains the crossing with $\widetilde D.$
		It follows that $Q'$ is divisible by $s_F^{a+\tau}$ as a polynomial on the affine logarithmic jet space.
		For a holomorphic arc $\gamma$ centered on $F$ and not contained in $D',$ let $p_F=\ord_0\gamma^*F.$
		The remaining polynomial of source weight $w$ has pole order at most $w,$ since each logarithmic derivative of order $h$ has pole order at most $h.$
		Hence $$\operatorname{pole}_0 Q'(j_k\gamma)\le\left(w-\left(a+\left\lceil\frac{2a}{d-3}\right\rceil\right)p_F\right)_+.$$
	\end{example}
	
		The hypersurfaces in the preceding higher-order examples are themselves Kobayashi hyperbolic by \cite[Theorem C]{DengDT}.
		A nonconstant entire curve in the blowup that is not contained in the boundary would give a nonconstant entire curve on the original hypersurface by composition with the blowup map.
		The comparison of pole bounds in the blowup example concerns local holomorphic arcs.
	
	\subsection{Sharpness and  comparison}
	
	The following product calculation tests the tangential part of the integer threshold and shows that its ceiling cannot be uniformly increased.
	Take smooth curves $C,B$ of genus at least two, $p\in B,$ and integers $1\le s<m.$
	Set $$X=C\times B,\ \ \ \ D=C\times\{p\},\ \ \ \ L=\operatorname{pr}_C^*K_C^s\otimes\operatorname{pr}_B^*(K_B(p))^{m-s}.$$
	Here, $L$ is ample, $D^2=0,$ $\mu=2g(C)-2,$ and $L\cdot D=s\mu.$
	The splitting $$\Omega_X^1(\log D)=\operatorname{pr}_C^*K_C\oplus\operatorname{pr}_B^*K_B(p)$$ gives a section with $s$ tangential factors and $m-s$ normal logarithmic factors.
	Its pole along a transverse arc with nonzero tangential derivative is exactly $m-s.$
	Replacing $K_C^s$ in $L$ by $K_C^s(-A),$ with $A$ effective and $0<\deg A<2g(C)-2,$ realizes nonintegral ratios with the same ceiling and the same pole outside $A.$
	Since the product has no nonconstant entire curves, this example tests the algebraic pole estimate.
	
	A full-weight pole at a crossing is also possible.
	For $X=(\PP^1)^2$ and $D$ the sum of $q\ge3$ vertical and $q$ horizontal fibers, $$\Omega_X^1(\log D)=\operatorname{pr}_1^*\mathscr O(q-2)\oplus\operatorname{pr}_2^*\mathscr O(q-2).$$
	Products using positive integer numbers $r_1,r_2$ of factors from the two summands give negatively ample-twisted sections of weight $m=r_1+r_2$ with pole order $m$ at a transverse crossing.
	The ample line bundle used for the twist is $$L=\operatorname{pr}_1^*\mathscr O_{\PP^1}\bigl(r_1(q-2)\bigr)\otimes\operatorname{pr}_2^*\mathscr O_{\PP^1}\bigl(r_2(q-2)\bigr).$$
	There, $s_i+s_j=m,$ so Theorem~\ref{prop:cross} gives exactly $m.$
	This example exhibits sharpness for the local pole bound, and the pair fails \eqref{eq:base} at every jet order.
	Indeed, fix a regular $k$-jet at a point of the first factor outside its boundary and set all logarithmic jet coordinates of the second factor to zero.
	These jets form a regular family parametrized by the entire second factor, including its boundary points.
	Evaluation of any global invariant logarithmic jet differential with twist $H^{-a}$ along that family gives a section of the restriction of $H^{-a}$ to the fiber of the first projection.
	This line bundle has negative degree, so the corresponding regular tower points outside $D$ are base points for every positive power in \eqref{eq:base}.
	Classical projection arguments already cover the counting inequalities of this example.
	
	Condition~\eqref{eq:base} is the sufficient hypothesis of \cite[Corollary 4.8]{BD}; constructions such as \cite[Corollary 4.5]{BD} establish it in suitable high-degree families.
	For a smooth irreducible boundary and each triple $k,m,a$ satisfying \eqref{eq:base}, Corollary~\ref{cor:smooth} gives a counting coefficient strictly below $m/a,$ independently of the power $q.$
	Comparison with other results for a special pair must also account for their section constructions and counting conventions.
	
	Restricted jet constructions require an additional argument.
	For example, in Case~2 of the proof in \cite[Theorem 1.9, Section 6.5]{HHMX}, when the initial global differential vanishes along the curve, the second differential is constructed on a divisor $Z\subset X_2(D).$
	In the cases where a global logarithmic jet differential has nonzero evaluation, Theorem~\ref{thm:mult} gives the subtraction $1/(d-3)$ from its weight-to-twist ratio for a smooth plane boundary of degree $d.$
	Take the inverse image of $Z$ under the quotient map from regular affine logarithmic $2$-jets to the regular-jet locus of $X_2(D),$ and let $\widetilde Z$ be its scheme-theoretic closure in $\mathcal J_2.$
	With the filtration induced on $\mathscr O_{\widetilde Z}$ by Lemma~\ref{lem:filtration}, its associated graded algebra is a quotient of the ambient associated graded algebra.
	Quotients need not preserve the maximal-slope upper bound or the absence of global sections, and the restricted locus can be singular along the pure normal locus.
	Therefore,  the present proof does not allow every coefficient $C_d$ in that work to be replaced by $C_d-1/(d-3).$
	A separate theorem for the restricted locus would be required.

\end{document}